\documentclass[oneside,english,american]{amsart}
\usepackage[T1]{fontenc}
\usepackage{textcomp}
\usepackage[utf8]{inputenc}
\usepackage{babel}
\usepackage{units}
\usepackage{amstext}
\usepackage{amsthm}
\usepackage{amssymb}
\usepackage{graphicx}
\usepackage{geometry}
\usepackage[pdfusetitle,
 bookmarks=true,bookmarksnumbered=false,bookmarksopen=false,
 breaklinks=false,pdfborder={0 0 1},backref=false,colorlinks=false]
 {hyperref}

\makeatletter
\numberwithin{equation}{section}
\numberwithin{figure}{section}
\theoremstyle{plain}
\newtheorem{thm}{\protect\theoremname}
\theoremstyle{definition}
\newtheorem{problem}[thm]{\protect\problemname}
\theoremstyle{remark}
\newtheorem{claim}[thm]{\protect\claimname}
\theoremstyle{definition}
\newtheorem{example}[thm]{\protect\examplename}
\theoremstyle{remark}
\newtheorem{rem}[thm]{\protect\remarkname}
\theoremstyle{definition}
\newtheorem{defn}[thm]{\protect\definitionname}
\theoremstyle{plain}
\newtheorem*{thm*}{\protect\theoremname}

\makeatother

\addto\captionsamerican{\renewcommand{\claimname}{Claim}}
\addto\captionsamerican{\renewcommand{\definitionname}{Definition}}
\addto\captionsamerican{\renewcommand{\examplename}{Example}}
\addto\captionsamerican{\renewcommand{\problemname}{Problem}}
\addto\captionsamerican{\renewcommand{\remarkname}{Remark}}
\addto\captionsamerican{\renewcommand{\theoremname}{Theorem}}
\addto\captionsenglish{\renewcommand{\claimname}{Claim}}
\addto\captionsenglish{\renewcommand{\definitionname}{Definition}}
\addto\captionsenglish{\renewcommand{\examplename}{Example}}
\addto\captionsenglish{\renewcommand{\problemname}{Problem}}
\addto\captionsenglish{\renewcommand{\remarkname}{Remark}}
\addto\captionsenglish{\renewcommand{\theoremname}{Theorem}}
\providecommand{\claimname}{Claim}
\providecommand{\definitionname}{Definition}
\providecommand{\examplename}{Example}
\providecommand{\problemname}{Problem}
\providecommand{\remarkname}{Remark}
\providecommand{\theoremname}{Theorem}

\begin{document}
\title{Cognitive Training for Language Models: \\
Towards General Capabilities via Cross-Entropy Games}
\date{March 24th, 2026}
\author{Cl\'{e}ment Hongler\textsuperscript{1,2}\quad Franck Gabriel\textsuperscript{3,1}\quad Valentin Hartmann\textsuperscript{1}\quad Arthur Renard\textsuperscript{1}\quad Andrew Emil\textsuperscript{1} \\[6pt]
\textsuperscript{1}Xent Labs\quad \textsuperscript{2}EPFL\quad \textsuperscript{3}Universit\'{e} Lyon 1}
\begin{abstract}
Defining a constructive process to build general capabilities for
language models in an automatic manner is considered an open problem
in artificial intelligence. Towards this, we consider the problem
of building a curriculum of tasks that grows a model via \emph{relevant
skill discovery}. 

We provide a concrete framework for this task, using a family of tasks
called Cross-Entropy Games, which we postulate is universal in a suitable
sense. We show that if it is possible to grow the curriculum for relevant
skill discovery by iterating a greedy optimization algorithm, then,
under natural assumptions, there is essentially only one \emph{meta-objective}
possible (up to a few hyper-parameters). We call the resulting process
\emph{cognitive training}. 

We postulate that, given sufficiently capable language models as players
and meta-samplers, cognitive training provides a principled way to
relevant skill discovery; and hence to the extent general capabilities
are achievable via greedy curriculum learning, cognitive training
would be a solution. 
\end{abstract}

\maketitle

\section{Artificial General Intelligence and Cognitive Training}\label{sec:artificial-general-intelligence-and-cognitive-training}

The last decades saw spectacular progress on several fronts for AI.
In particular we saw:
\begin{itemize}
\item Models that could generalize beyond their training dataset.
\item Reinforcement learning that could reach super-human performance on
specific tasks.
\item Pre-trained Large Language Models (LLMs) that could learn to perform
certain tasks with in-context learning, and numerous more specific
tasks, given a suitable training environment. 
\end{itemize}
In spite of these achievements, there is a growing consensus that
a number of ideas are lacking to go towards ``true Artificial General
Intelligence (AGI)''. In this paper, we start with a short summary
of what we expect from AGI starting from seminal works \cite{turing,legg-hutter}
and propose to study an (apparently simpler) problem, which is that
of (open-ended) relevant skill discovery. 

To study this problem, we propose a training framework (based on an
elaboration of the framework introduced in \cite{hongler-emil}) which
we call \emph{cognitive training}: the idea is to grow a curriculum
of games in a specific space of tasks, called \emph{Cross-Entropy
(Xent) Games}, endowed with a suitable notion of \emph{transfer value}:
how learning one game teaches one to play another. We postulate that
the learning achieved by any curriculum of tasks can be approximated
by a curriculum of Xent Games. 

We consider the problem of growing a curriculum of Xent Games in a
greedy manner, by optimizing a meta-objective $\mathcal{O}$. We show
that (perhaps surprisingly), given a few natural assumptions, a consistent
meta-objective must take a very constrained form: from these, we derive
an explicit formula for $\mathcal{O}$, that balances sparsity (inverse
code length), quality (internal consistency), diversity (novelty),
and external relevance (benchmark performance).

This leads us to \emph{cognitive training}. The idea is to start with
one (or more) base auto-regressive model $\mathcal{M}$, find a suitable
space of Xent Games $\mathcal{G}$. At every step of the curriculum,
a meta-sampler $\mathcal{M}_{\mathcal{S}}$ then generates games on
$\mathcal{G}$ to optimize $\mathcal{O}$, to grow a curriculum that
brings relevant skill discovery to $\mathcal{M}$. 

Our reasoning thus leads (from a number of assumptions) to the following
idea: \emph{if it is possible to greedily generate a curriculum of
tasks bringing general capabilities to a model, then it is possible
to replicate this with a greedily built curriculum of Xent Games,
and then the meta-objective formula must take a very specific form. }

\subsection{Definition of Artificial General Intelligence}\label{subsec:definition-of-agi}

The first notable attempt at a measure of artificial intelligence
is probably Turing's imitation game \cite{turing}, leading to the
famous Turing test. Across decades, various definitions of intelligence
in the context of AI were provided, leading to a notable and influential
synthesis in \cite{legg-hutter}:

\emph{Intelligence measures an agent’s ability to achieve goals in
a wide range of environments.}
\begin{flushright}
S. Legg and M. Hutter
\par\end{flushright}

A framework built at the intersection of Occam's Razor (in the form
of Algorithmic Probability) and of Reinforcement Learning \cite{legg-hutter-2}
highlighted notably a theoretically well-posed, but uncomputable,
approach to the so-called Universal AGI problem: the AIXI paradigm
\cite{hutter}, blending a universal prior on world models (based
on Kolmogorov's complexity) and objective maximization. 

Arguably, the part that remains a little vague is associated with
the ``wide range'' in the above definition: is there a good measure
of generality? Are humans general according to any such definition?
Is there a process through which this generality is discovered?

Despite the impossibility to apply this framework directly (due to
its uncomputability) and the lack of answers about generality, AIXI
gives a direction for what we would expect theoretically from an AI
going towards AGI: models should interact with an environment, perform
tasks within them, and be rewarded. Two decades of progress in the
field of Reinforcement Learning (RL) followed. In the context of language
models, this was discussed for quite some time \cite{mikolov-joulin-baroni}.
However, as discussed in the next subsection, the progress has only
recently accelerated thanks to the rise in prominence of pre-trained
large language models. 

\subsection{Cognitive Training: Goals}\label{subsec:cognitive-training-goals}

The current problem with language models is not that they cannot learn
tasks given a suitable RL environment, but that they have a hard time
performing (or quickly learning to perform) tasks for which they have
not been trained already. Given a new task, constructing a specifically
relevant RL environment and post-training a model on it can be very
difficult or impossible. 

The question that we aim to address can be understood as building
extended pre-training using \emph{games} (tasks/environments with
an emphasis on certain types of rewards):
\begin{problem}
\label{prob:find-small-family-of-games}Can we find a curriculum of
compact text games $\left(G_{k}\right)_{k\geq0}$ such that a model
$\mathcal{M}$ trained on it is maximally ready for new/unseen tasks
given some data, architecture, and compute constraints?
\end{problem}

Mathematically, the above problem seems to assume a prior on new tasks,
and is hence somewhat ill-posed (we have no way to set this prior).
Like for supervised learning (where e.g. computer vision problems
could be solved without a clear prior on the set of images corresponding
to a given label), it is reasonable to expect that one can achieve
good performance with the right methods, without expliciting (or sampling
from) such a prior. By Occam's razor, if we find a curriculum of games
(at finite $k$) that yields improved performance on an external metric
$E$ (e.g. an aggregation of benchmarks), then it is reasonable to
expect that this curriculum will be \emph{generally} valuable: a model
trained on such games is likely (in an informal sense) to generalize
well to tasks beyond the specific tasks of $E$.

Pre-training is in fact the most emblematic example of a compact training
game: predicting the next token, a most universal and abstract task,
leads us quite far in terms of general capabilities. Given a dataset,
this is, however, a game that can only be played once (and with diminishing
returns) leading to the question: what are, besides pre-training,
the most useful training games?

Interestingly, we can perform a substantial reduction to the above
question by asking the following, which seems perhaps a simpler problem
a priori: 
\begin{problem}
\label{prob:find-open-ended-skill-discovery}Can we find a curriculum
of compact games $\left(G_{k}\right)_{k\geq0}$ such that an LLM $\mathcal{M}$
will keep discovering relevant new skills throughout training, while
maintaining existing capabilities?
\end{problem}

The following seems intuitive, as any solution to Problem \ref{prob:find-small-family-of-games}
would (at least eventually) involve the discovery of new relevant
skills (which are needed to deal with new tasks):
\begin{claim}
\label{claim:requirement-for-agi}Any convincing solution to Problem
\ref{prob:find-small-family-of-games} should also yield a solution
to Problem \ref{prob:find-open-ended-skill-discovery}. 

According to Claim \ref{claim:requirement-for-agi}, finding a curriculum
as described in Problem \ref{prob:find-open-ended-skill-discovery}
is a prerequisite for training a model that is maximally prepared
for new tasks. Deriving a process for finding such a curriculum is
the goal of this work.

To make the above problem concrete, we propose a space of games, called
\emph{Cross-Entropy Games (Xent Games),} endowed with an appropriate
structure of transfer value; separately, it is postulated that any
curriculum of tasks can be approximated by a curriculum of Xent Games
(Section \ref{subsec:xent-game-space-properties}); this leads to
a specific framework for the above definition (Definition \ref{def:relevant-skill-discovery}). 
\end{claim}

We then focus on \emph{greedy curriculum building} (i.e. where each
game is picked by optimizing a function of the previous games). This
results in the central (and perhaps most surprising) claim of this
paper (Sections \ref{subsec:meta-objective-formula}--\ref{subsec:meta-objective-derivation}): 
\begin{claim}
If we can solve Problem \ref{prob:find-open-ended-skill-discovery}
via a greedy method involving at every step the optimization of meta-objective
$\mathcal{O}$, then from reasonable principles an explicit (and unique)
formula for $\mathcal{O}$ can be found.
\end{claim}

This leads us to the concept of \emph{cognitive training}: the process
of performing the meta-optimization for $\mathcal{O}$. In the next
subsection, we first lay out a framework to optimize cognitive training
for LLMs, detailed in the rest of the paper.

\subsection{Cognitive Training: Implementation}\label{subsec:cognitive-training-implementation}

The question raised in Problem \ref{prob:find-small-family-of-games}
above suggests formulating a meta-optimization problem: that of building
a curriculum of games that turn out to be maximally useful to a (language)
model $\mathcal{M}$ learning to play them (for a good definition
of ``usefulness''), on which one trains a model sequentially. 

To implement this program, we need two ingredients:
\begin{itemize}
\item a space of games $\mathcal{G}$ on which the model $\mathcal{M}$
will be trained,
\item a meta-objective $\mathcal{O}$ defined on games $\mathcal{G}$, used
to select which games to add to the curriculum.
\end{itemize}
Given these two ingredients, $\mathcal{O}$ is optimized at each step
by a meta-sampler $\mathcal{M}_{\mathcal{S}}$, which outputs games
(written in a specific language). The next game added to the curriculum
to train model $\mathcal{M}$ is the one that maximizes the meta-objective
$\mathcal{O}$. In the rest of this subsection, we briefly discuss
principles for the design of $\mathcal{G}$ and $\mathcal{O}$; details
are presented in Sections \ref{sec:streamlined-xent-games} and \ref{sec:cognitive-training},
respectively.

\subsubsection{Game Space Design Elements}\label{subsec:game-space-design-elements}

The desiderata for the space of games $\mathcal{G}$ aimed at developing
the cognitive capabilities of a model $\mathcal{M}$ are numerous: 
\begin{itemize}
\item The games should be amenable to fast training.
\item It should be easy to generate a wide diversity of games in $\mathcal{G}$
with fairly concise code. 
\item The space should allow for robust transfer exploration: if a certain
curriculum (based on certain games not necessarily in $\mathcal{G}$)
teaches $\mathcal{M}$ some skills, the same can be achieved by a
curriculum of games in $\mathcal{G}$.
\end{itemize}
In Section \ref{sec:implicit-knowledge} below, we argue that a fertile
ground for building the games in $\mathcal{G}$ is to rely on the
\emph{implicit knowledge} of LLMs; this leads to the proposal to use
Cross-Entropy Games (Xent Games) as a means from which to build $\mathcal{G}$. 

\subsubsection{Meta-Game Design Elements}\label{subsec:meta-game-design-elements}

Our main contribution is the principled derivation in Section \ref{sec:cognitive-training}
of a meta-algorithm, called \emph{cognitive training}, to build a
useful curriculum of games $\left(G_{k}\right)_{k\geq0}$ in $\mathcal{G}$
to develop the cognitive capabilities of a model $\mathcal{M}$, building
a sequence of models where $\mathcal{M}_{k+1}$ is obtained by training
$\mathcal{M}_{k}$ on $G_{k}$. The \emph{cognitive training} algorithm
is a greedy algorithm based on a \emph{meta-objective }$\mathcal{O}$.
Central to the definition of $\mathcal{O}$ is the concept of \emph{transfer
value} (Section \ref{subsec:transfer-value}), which estimates the
relevance of learning a game to play another. 

At each step $k\geq1$, the choice of a new game $H$ should balance:
\begin{itemize}
\item quality/relevance $q$: how much $H$ improves performance on the
old games $G_{<k}:=G_{1},\dots,G_{k-1}$;
\item diversity/novelty $d$: how much $H$ brings skills that are not yet
captured by the games $G_{<k}$;
\item external benchmark performance $b$: how much $H$ improves performance
on external metrics;
\item description length $l$: the raw code length of $H$. 
\end{itemize}
In Section \ref{sec:cognitive-training}, we derive explicit principled
definitions of each of the previous quantities, and derive a formulation
for the meta-objective $\mathcal{O}$, namely 
\[
\mathcal{O}=\frac{qd+bp}{l},
\]
where $p$ is a \emph{pressure} factor, balancing the relative importance
of internal and external performance. Note that despite involving
a priori a number of terms that scale linearly in $k$ (and similarly
for memory), the meta-objective $\mathcal{O}$ satisfies some good
scalability properties (see Section \ref{subsec:scalability-of-meta-objective-computation}). 

\subsubsection{Scope and Contributions}

The key contributions of this paper are the following:
\begin{itemize}
\item A treatment of Problem \ref{prob:find-open-ended-skill-discovery},
i.e. relevant skill discovery, with the idea that if general capabilities
are achieved by curriculum learning, they automatically imply relevant
skill discovery.
\item A formal definition, building upon \cite{hongler-emil}, of a space
of tasks called \emph{(Streamlined) Cross-Entropy (Xent) Games}, endowed
with a transfer value structure.
\item A hypothesis that curricula made of (Streamlined) Xent Games can approximate
the learning of any curriculum; and hence that relevant skill discovery
can be understood in terms of Xent Games and their transfer value
structure. 
\item A formalization of greedy curriculum learning for Xent Games in terms
of a meta-objective optimization.
\item A derivation of an explicit formula for any meta-objective satisfying
suitable consistency assumptions.
\end{itemize}

\subsubsection{Structure of the Paper}
\begin{itemize}
\item In Section \ref{sec:implicit-knowledge}, we outline the key idea
of implicit knowledge for an LLM: this suggests that much more can
be done with an LLM than simple sampling.
\item In Section \ref{sec:streamlined-xent-games}, we present the space
of games upon which cognitive training is based: the (Streamlined)
Xent Games, derived from the implicit knowledge of LLMs.
\item In Section \ref{sec:cognitive-training}, we derive a natural framework
for xent-game-based cognitive training for LLMs.
\item In the same Section \ref{sec:cognitive-training}, we discuss several
challenges and open questions associated with cognitive training. 
\end{itemize}

\section{Implicit Knowledge}\label{sec:implicit-knowledge}

If the models that we train are supposed to be useful for general
tasks relevant to some world (e.g. human society), the games that
they are trained on need to have a connection to that world. For language
models, we posit that this connection is provided by the \emph{implicit
knowledge} of Large Language Models (LLMs) that have been trained
on vast corpora of text data.

LLMs are generative models: given a context $c$, they define a conditional
probability measure $\mathbb{P}_{\mathcal{M}}(\cdot\mid c)$ on sequences
of tokens. For a given model $\mathcal{M}$, the likelihood is directly
accessible: for any token sequence $x=x_{1},\ldots,x_{T}$, the quantity
$\log\mathbb{P}_{\mathcal{M}}(x\mid c)=\sum_{t=1}^{T}\log\mathbb{P}_{\mathcal{M}}(x_{t}\mid c,x_{<t})$
is readily available from $\mathcal{M}$. As a consequence, the information
we can extract from $\mathcal{M}$ depends on how we use it: either
as a generator, or as the basis for an algorithm relying on $\mathbb{P}_{\mathcal{M}}$
for e.g. perplexity comparison, search, and optimization in the space
of token sequences. 

In this section, we discuss the difference between: 
\begin{itemize}
\item Explicit knowledge $\mathcal{E}_{\mathcal{M}}$: the information that
can be quickly and reliably extracted from $\mathcal{M}$ by sampling
answers from well-chosen prompts. 
\item Implicit knowledge $\mathcal{I}_{\mathcal{M}}$: the information that
can be extracted by arbitrary algorithms with access to $\mathbb{P}_{\mathcal{M}}$.
\end{itemize}
Since generation is itself a (randomized) algorithm relying on $\mathbb{P}_{\mathcal{M}}$,
we have $\mathcal{E}_{\mathcal{M}}\subset\mathcal{I}_{\mathcal{M}}$.
In the following, we will see that there is a gap between the two
kinds of knowledge: focusing on continuations allows one to explore
only a tiny subset of observables (functions) of the measure $\mathbb{P}_{\mathcal{M}}$,
while direct access to the model's measure allows for different procedures
(search, contrastive experiments, ...) that further reveal the cognitive
capabilities of the LLM. In Section 3, we exploit the implicit/explicit
gap to define the Xent Games, designed to use implicit capabilities
to improve the explicit knowledge of LLMs. 

\subsection{Explicit Knowledge}\label{subsec:explicit-knowledge}

The \emph{explicit knowledge} of an LLM is defined as follows: the
set of questions and tasks for which, given a fixed computation budget,
a model $\mathcal{M}$, typically fine-tuned for instruction following,
can provide a correct solution with sufficiently high probability
by generating a continuation. This knowledge is essentially what a
standard benchmark (e.g. one based on the samples generated by the
LLM) would measure. 

The notion of explicit knowledge has some important subtleties: 
\begin{enumerate}
\item The explicit knowledge for raw pre-trained models (not instruction-tuned)
is ill-defined or poorly representative of the underlying knowledge
of the LLM: for such models, the main objective is not to answer or
to follow the output specifications. \footnote{With well-chosen prompts, one can steer LLMs towards the desired behavior,
but this will not replace a good fine-tuning.}
\item Explicit knowledge depends on the decoding policy: it measures both
the model distribution and an external choice of how we generate from
it (greedy/sampling, temperature, ...).
\item The answers of an LLM can be stochastic, so we may allow for the sampling
of a few answers before an algorithm produces the final output from
these samples. 
\end{enumerate}
Finally, fine-tuning can modify part of the explicit knowledge: for
instance, after alignment training, a model may refuse to answer a
question its base version could have answered. This motivates the
following perspective for understanding LLMs' capabilities: instead
of asking the model what it knows, we can ask what can be extracted
from the model probabilities. 

\subsection{Implicit Knowledge}\label{subsec:implicit-knowledge}

We define the \emph{implicit knowledge} as follows: the set of questions/tasks
that can be answered algorithmically (given a fixed computation budget)
given access to the model measure $\mathbb{P}_{\mathcal{M}}$. More
specifically, we assume that we have access to the log-likelihoods
of continuations $\log\mathbb{P}_{\mathcal{M}}(x\mid c)$. The notion
of implicit knowledge is obviously appealing as it is the set of information
that an LLM already knows (in one form or another) from its pre-training. 

The implicit knowledge differs from the explicit knowledge in at least
three important ways. First, it is well defined for raw pre-trained
models, as it does not rely on a question-answering setting. Second,
it does not depend on the decoding policy as it depends solely on
$\mathbb{P}_{\mathcal{M}}$. Third, it is, by its nature, related
to a deterministic computational problem, rather than the consequence
of a stochastic generator.

The simplest but most fundamental implicit knowledge is the ability
to provide the log-likelihood of a sentence or continuation. For any
$x=x_{1}\ldots x_{T}$, 
\[
\log\mathbb{P}_{\mathcal{M}}\left(x_{1}\ldots x_{T}\mid c\right)=\sum_{t=1}^{T}\log\mathbb{P}_{\mathcal{M}}\left(x_{t}\mid c,x_{<t}\right),
\]
can be extracted from the model's measure, and is therefore part of
the implicit knowledge. Yet, this quantity is not part of the explicit
model: a model generally cannot explicitly report its correct value. 

Of course, by sampling more and more continuations, we should theoretically
be able to recover $\log\mathbb{P}_{\mathcal{M}}(x\mid c)$. Yet,
the number of samples required is prohibitive. Indeed, for a fixed
continuation $x$, observing $x$ even once requires an order of $\nicefrac{1}{\mathbb{P}_{\mathcal{M}}(x\mid c)}$
samples. Since $\mathbb{P}_{\mathcal{M}}(x\mid c)$ typically decreases
exponentially with the length of $x$, it is impossible to recover
the log-likelihood from samples within a reasonable computing time.
This simple argument shows that having direct access to the log-likelihood
is qualitatively different from having access only to the generation
process. 

\subsection{Examples of Implicit Knowledge of Pre-Trained Models}\label{subsec:examples-of-implicit-knowledge-of-pre-trained-llms}

Direct access to the log-likelihood enables procedures that are out
of reach for explicit generation. As explained above, this is important
for raw pre-trained models: even when a model does not reliably follow
a question/answer format, its conditional measure can still be evaluated
through log-likelihoods. In this section, we provide a non-exhaustive
family of procedures for extracting implicit knowledge. 

One way to do so is to consider the difference between log-likelihoods
in which some candidate sequences $x_{1}$ or $x_{2}$ are inserted
into possibly different contexts $c_{1}$ and $c_{2}$. This quantity
measures how much more the model supports $x_{1}$ in the hypothetical
world specified by $c_{1}$ than $x_{2}$ in the hypothetical world
specified by $c_{2}$. By choosing contexts appropriately and determining
how to insert the candidate sequences into them, this likelihood difference
can be used among other things to extract beliefs from a model, quantify
the usefulness of some information, or define some contrastive objectives.
In practice, contexts are built from a base context to which pre-prompts
or sandwich prompts are added before and after it; we omit these details
here, but they matter in concrete implementations. 

\subsubsection{Likelihood difference for a fixed statement}\label{subsec:likelihood-ratio-statement}

We can compare the likelihood of the same statement in two different
contexts. 
\begin{itemize}
\item \textbf{Truth-false difference.} For a statement $s$, we can consider
\begin{eqnarray*}
\Delta_{\mathrm{T/F}}(s\mid c) & = & \log\mathbb{P}_{\mathcal{M}}(s\mid c,"\text{The following statement is true}:")\\
 &  & -\log\mathbb{P}_{\mathcal{M}}(s\mid c,"\text{The following statement is false}:")
\end{eqnarray*}
A simple classifier predicts ``true'' whenever $\Delta_{\mathrm{T/F}}(s\mid c)>0$.
To reduce dependence on the exact prompt, we can compute the average
of $\Delta_{\mathrm{T/F}}$ over paraphrases of the prefix and/or
of the statement $s$. Closely related, in \cite{Language Models (Mostly) Know What They Know},
to study the calibration of models, the authors present a proposed
answer to a model and ask it whether it is true or false, then read
off the log-likelihood of the two options. 
\item \textbf{Multiple-choice selection for explanations.} More generally,
given candidate explanations $c_{1},\ldots,c_{n}$, we can select
$\arg\max_{i\in\{1,\ldots,n\}}\log\mathbb{P}_{\mathcal{M}}(s\mid c_{i})$.
This allows one to select the most plausible explanation for $s$. 
\end{itemize}

\subsubsection{Counterfactual thinking and surprise reduction}\label{subsec:Counterfactual-thinking}

A variant of the multiple-choice selection of explanations is to quantify
the value of some information. Let $q$ be a problem, $x$ a reference
solution, and $h$ an additional piece of information. The difference
\[
\Delta_{\mathrm{info}}(h;q\to x)=\log\mathbb{P}_{\mathcal{M}}(x\mid q,h)-\log\mathbb{P}_{\mathcal{M}}(x\mid q)
\]
measures how much $h$ reduces the model's surprise about $x$. This
makes it possible to quantify the usefulness of a piece of information
without requiring the model to explicitly explain why (which does
not provide a reliable, quantitative value). 

\subsubsection{Contrastive objective across models or checkpoints }\label{subsec:Contrastive-objective}

Another variant of log-likelihood comparison is to consider two different
models or checkpoints. Given two models $\mathcal{M}_{1}$ and $\mathcal{M}_{2}$,
the quantity
\[
\Delta_{\mathrm{ctr}}(x\mid c)=\log\mathbb{P}_{\mathcal{M}_{1}}(x\mid c)-\log\mathbb{P}_{\mathcal{M}_{2}}(x\mid c)
\]
quantifies whether $x$ is more plausible for $\mathcal{M}_{1}$ than
for $\mathcal{M}_{2}$. When $\mathcal{M}_{1}$ is stronger, for example
if it is a larger model, and $\mathcal{M}_{2}$ is weaker, this quantifies
the continuations that are more interesting in the sense that a clever
model can think about them, but a simpler model cannot. This is closely
related to contrastive sampling methods (see e.g. \cite{contrastive}). 

\subsubsection{Verification vs construction}\label{subsec:verification-vs-construction}

Many tasks have solutions that are easier to verify than to construct.\footnote{In the context of self-evaluation, the authors of \cite{Language Models (Mostly) Know What They Know}
note that ``verification improves faster than generation quality
in this context''. } A well-trained model can assign a higher likelihood to correct solutions
than to incorrect ones, while still failing to construct a correct
solution by direct generation. 

Access to the log-likelihood transforms the verification task into
a continuous signal: rather than a discrete ``correct/incorrect''
reward, $\log\mathbb{P}_{\mathcal{M}}$ can be used to construct a
continuous objective that can guide search and optimization. 

In addition, it is worth noting that some continuity in reward signals
can be provided at the length level: rewards on incomplete responses
can be computed as well, and help steer the trained models towards
high-reward responses. 

This perspective can be useful to learn to solve problems that can
even be a priori combinatorial in nature (by providing some suitable
continuous relaxations of such problems using the LLM cross-entropy
function to provide a ``soft'' enforcement of the constraints).
The gap between verification and construction is one of the important
gaps that can produce new data and games which can later be distilled
into the explicit knowledge. 

\subsubsection{Constrained optimization in the space of sequences and inverse prompting}

More generally, the log-likelihood (together with suitable contexts)
can be used to define various losses on token sequences. A fundamental
optimization problem is the prompting problem, which consists of finding
the optimal set of tokens that satisfies some fixed constraints (format,
length, vocabulary restrictions, ... ) and maximizes the log-likelihood
when inserted into a fixed context. This family includes some meaningful
inverse problems: 
\begin{example}[Inverse prompting and summarization]
\label{exa:inverse-prompting}Given a fixed text $x$, we can search
for a short sequence of tokens $t$ such that $x$ becomes highly
likely when it follows $t$, i.e. 
\[
t^{*}\in\arg\max_{t\in\mathcal{C}}\log\mathbb{P}_{\mathcal{M}}(x\mid t),
\]
where $\mathcal{C}$ encodes the constraints. With appropriate prompting
between, before or after $t$ and $x$, this can be interpreted as
a form of summary: $t$ summarizes $x$ because it reduces the surprise
of $x$ for the model. This type of implicit knowledge is a direct
building block of the subsequent Xent Games.
\end{example}

\begin{example}[Common explanations]
\label{exa:common-explanations}A variant of the previous example
can be obtained when given multiple texts $x_{1},\ldots,x_{n}$. We
can search for a sequence of tokens $t$ that makes them jointly more
likely, but at the same time relatively unexpected from each of them
individually 
\[
t^{*}\in\mathrm{arg}\max_{t\in\mathcal{C}}\log\mathbb{P}_{\mathcal{M}}(x_{1},\ldots,x_{n}\mid t)-\alpha\sum_{i=1}^{n}\log\mathbb{P}_{\mathcal{M}}(t\mid x_{i}).
\]
With appropriate prompting, $t$ can be understood as an interesting
common feature about the texts.
\end{example}

\subsubsection{Anomaly detection}

If a small part of a text has a high cross-entropy loss, it probably
deserves attention: this could signal an anomaly, an adversarial or
injected segment, or simply a rare but meaningful information in the
text that the model cannot well predict. This intuition has been used
in practice for e.g. adversarial prompt detection \cite{token-level-adversarial-prompt-detection}
and outlier spatial trajectories \cite{mbuya-pfosser-anastosopoulos}. 

\subsection{Explicit vs Implicit Knowledge}

As discussed, the explicit knowledge of an LLM is a subset of its
implicit knowledge, since sampling tokens from the model is an algorithm
on its probability measure. It is a strict subset, because none of
the probabilities in the examples in \ref{subsec:examples-of-implicit-knowledge-of-pre-trained-llms}
can be computed with a (reasonable) amount of sampling. While some
questions like the truth value of a statement can also be answered
by generating tokens, having access to the probability allows us to
determine the certainty of the model. This access through the model's
probability measure furthermore does not require fine-tuning for instruction-following
and is thus well-defined for pre-trained models. 

Due to its vastness, generality and easy accessibility, we use the
implicit knowledge of LLMs as the basis of Xent Games, introduced
in the next section. These games will make up the training curriculum
of an agent designed for general capabilities.

\section{Streamlined Xent Games}\label{sec:streamlined-xent-games}

Building on the idea of implicit knowledge outlined in Section \ref{sec:implicit-knowledge}
above, we briefly outline the ideas behind Xent (Cross-Entropy) Games.
Rather than fully specifying the language used to write Xent Games
(see \cite{hongler-emil} for details), we focus on a core subset
of them that we call \emph{Streamlined Xent Games (S-XGs)}, and we
provide illustrative examples. 

S-XGs form a subspace of the Original Xent Games (O-XGs) introduced
in \cite{hongler-emil}, which focuses on benchmarking. By contrast,
S-XGs are designed to make training more streamlined (i.e. more parallelizable,
easier to optimize, and with extra differentiable structures \cite{hongler-emil-renard-gabriel,hongler-asani-gabriel-emil-renard})
and therefore suitable for training models to play them. At the same
time, it is reasonable to expect that the set of skills developed
by S-XGs is the same as that developed by the space of O-XGs (see
Section \ref{subsec:web-of-games} below).

\subsection{Definition and Runtime}\label{subsec:definition}

Informally, Xent Games are text-based games played by some main players
$\mathcal{M}$, with one or more LLMs used as ``world models'' (i.e.
providing the environment dynamics and rewards). The whole space of
games is endowed with some meta-data, corresponding to the specification
of the models involved, which are each attached to a \emph{model name}
(i.e. they link variable names used in the game codes to concrete
model checkpoints). In the simplest variant (on which we focus here),
one of these is the main player model, which is the one under cognitive
training; the other models are frozen and they play the roles of judge,
data stream models, or NPCs (opponents / cooperators). 

The game state consists of a collection of string registers that are
updated according to basic string operations (Section \ref{subsec:token-string-space}),
by using inputs from data streams, and by writing players' outputs. 

The players receive as input strings read from the registers and produce
output strings that are written back into the registers. In addition,
they receive rewards based on signed cross-entropies of token strings
stored in the registers (Section \ref{subsec:xents}). 

\subsubsection{Token String Space}\label{subsec:token-string-space}

As text-based games, Xent Games run on a space of strings, endowed
with a small set of basic string operations. In the framework of O-XGs,
the allowed operations are: `cat' (i.e. string concatenations at the
character level) and `cut' operations (i.e. splits between what comes
before a first occurrence of a string and what comes after). In S-XGs,
the situation is simpler: 
\begin{itemize}
\item Strings each have a token length that can change over the course of
a game.
\item Concatenations and cuts are replaced by copy operations made at the
token level (rather than the character level) that copy part of a
token string into another. The append operation extends the length
of a string up to the maximal length, the cut operation moves tokens
from one string to another until the second string's length is reached.
\end{itemize}

\subsubsection{Xents and Xent Sums}\label{subsec:xents}

At the heart of Xent Games are cross-entropies of strings $\mathrm{xent}_{\mathcal{J}}$
computed by a (judge) model $\mathcal{J}$. If $x=x_{1},\ldots,x_{m}$
and $y=y_{1},\ldots,y_{m}$ are token strings (read from a token string
register), we define 
\[
\mathrm{xent}_{\mathcal{J}}\left(x|y\right)=-\log\mathbb{P}_{\mathcal{J}}\left(x|y\right)=-\sum_{i=1}^{n}\log\mathbb{P}_{\mathcal{J}}\left(x_{i}|y,x_{<i}\right).
\]
Informally, $\mathrm{xent}_{\mathcal{J}}\left(x|y\right)$ measures
how ``surprised'' the autoregressive model $\mathcal{J}$ is to
see $x$ after seeing $y$. We denote also $\mathrm{xent}_{\mathcal{J}}(x)=\mathrm{xent}_{\mathcal{J}}(x|\emptyset)$
where $\emptyset$ is the empty string. 

The (judge) model $\mathcal{J}$ can be the main player model $\mathcal{M}$
itself (as in, e.g., pre-training games; see below), or a fixed pre-trained
or instruct model. Depending on whether $\mathcal{J}=\mathcal{M}$
or $\mathcal{J}\neq\mathcal{M}$, the goal may be to improve the model
used for the cross-entropy computation, or simply to optimize with
respect to this cross-entropy defined by a fixed $\mathcal{J}$. 
\begin{rem}
In practice, for S-XGs, it can be useful to clip the $\mathbb{P}_{\mathcal{J}}\left(x_{i}|y,x_{<i}\right)$
from below by a small amount, which clips the per-token loss $-\log\mathbb{P}_{\mathcal{J}}\left(x_{i}|y,x_{<i}\right)$
from above. This prevents players from exploiting very low probabilities
of models to achieve high rewards: these log-probabilities are typically
noisy and not meaningful (see the well-posedness criterion in \cite{hongler-emil}). 
\end{rem}

For a family of token strings $\{(x_{j},y_{j})\}_{j}$ and a family
of models $(\mathcal{J}_{j})_{j}$ a signed \emph{xent sum} is an
expression of the form 
\[
\sum_{j}\sigma_{j}\mathrm{xent}_{\mathcal{J}_{j}}\left(x_{j}|y_{j}\right),\qquad\sigma_{j}\in\left\{ \pm1\right\} 
\]
Xent sums can be used in two ways:
\begin{itemize}
\item as rewards given to the players; 
\item as skewed rewards (i.e. soft constraints, see next Section \ref{subsec:sxg-def-and-basic-runtime-instruction})
given to the players. 
\end{itemize}

\subsubsection{S-XG: Definition and Basic Runtime Instruction}\label{subsec:sxg-def-and-basic-runtime-instruction}

An S-XG consists of a sequence of moves of four types:
\begin{itemize}
\item \emph{assign}: modify a string register using the token-string operations
described in Section \ref{subsec:token-string-space}. 
\item \emph{elicit}: request from a model a string of length $n$ tokens. 
\item \emph{reward}: reward a player based on a signed xent sum (evaluated
on token-string registers).
\item \emph{ensure}: provide a smooth constraint on a player by imposing
a soft positivity constraint $S\geq0$ on a signed xent sum $S$ (evaluated
on token-string registers). Concretely, for a fixed $\lambda>1$,
the model receives reward $S/\lambda$ if $S\geq0$ or $\lambda S$
if $S<0$ (see Remark \ref{rem:ensure} below). 
\end{itemize}
\begin{rem}
\label{rem:ensure}The S-XG \emph{ensure} constraints follow the same
``adversarial multiplier'' logic as in O-XG \cite{hongler-emil},
but with a fixed finite multiplier $\lambda>1$ in case of violation
(and a small nonzero multiplier $1/\lambda$ in case of fulfillment).
Natural examples of \emph{ensure} constraints are true/false statements
based on the implicit knowledge (as in Section \ref{subsec:likelihood-ratio-statement}).
Another example of an \emph{ensure} constraint is based on the signed
xent sum 
\[
\mathrm{xent}_{\mathcal{J}}(s2|s1+\text{"comes after"})-\mathrm{xent}_{\mathcal{J}}(s1|s2+\text{"comes after"}),
\]
which softly enforces that $s2$ is more likely to follow $s1$ than
the reverse. 
\end{rem}

\subsubsection{Original Xent Game Language (O-XGL)}\label{subsec:oxgl}

The simple principles behind O-XGs make them naturally suited to being
expressed in a domain-specific language, which we call \emph{O-XGL}
(Original Xent Game Language). The design is close in spirit to an
assembly language: each line of code corresponds to one instruction
of the game logic. Any program made of valid lines of code is valid
(there are no conditional jumps), which allows for e.g. a batched
execution for training. Each line corresponds to a basic instruction
of the types listed in Section \ref{subsec:sxg-def-and-basic-runtime-instruction}
above. In the O-XGL specification \cite{hongler-emil}, a number of
instructions made to simplify the writing and reading of Xent Games
by humans are added; however, they do not alter the space of games
being considered, and for simplicity, we will omit those here. 

The inputs to the \emph{reward} and \emph{ensure} statements are signed
xent sums. In O-XGL, these instructions take three ingredients as
inputs: the judge, the target string, and the (possibly empty) prefix
strings. 

\subsubsection{Streamlined Xent Game Language (S-XGL)}\label{subsec:sxgl}

As explained above, for the purpose of cognitive training, we rely
on streamlined Xent Games (S-XGs). This motivates the S-XGL (Streamlined
Xent Game Language) specification. S-XGL is a minimalistic, assembly-like
language, which only consists of two basic binary operators, $<<$
and $>>$, acting on two different types of objects: 
\begin{itemize}
\item Models (main player, judges, data streams, NPCs). Each model has a
context register and a score register, used to accumulate xent-based
scores and to implement the ensure moves.
\item Token strings with individual lengths and a global fixed maximal length,
implementing the token-space operations described above (see Section
\ref{subsec:token-string-space}).
\end{itemize}
For example, given a model $\mathcal{M}$ and a string register $\mathcal{S}$: 
\begin{itemize}
\item $\mathcal{S}<<\mathcal{M}$ samples from $\mathcal{M}$ and appends
the resulting tokens to $\mathcal{S}$ until the length of $\mathcal{S}$
has doubled or the maximal length is reached,
\item $\mathcal{M}<<\mathcal{S}$ appends $\mathcal{S}$ to the context
of $\mathcal{M}$,
\item $\mathcal{S}>>\mathcal{M}$ adds $\mathrm{xent}_{\mathcal{M}}\left(\mathcal{S}\right)$
to the xent accumulator of $\mathcal{M}$, and$\mathcal{M}>>\mathcal{S}$
subtracts it from the xent accumulator,
\item $\mathcal{M}_{\ell}<<\mathcal{M}_{r}$ rewards $\mathcal{M}_{\ell}$
with the value in the xent accumulator of $\mathcal{M}_{r}$, and
clears that xent accumulator. 
\end{itemize}
Since the two binary operators can act on an ordered pair of objects
(of two possible types), this results in $2\times2\times2=8$ types
of operations. Together with the empty line operation, which clears
all strings, context registers and score registers, this makes $8+1$
operations, which constitute all of S-XGL instructions. All other
types of lines are valid code (with no execution value), which can
be loaded into token strings using S-XGL instructions. A detailed
description of S-XGL is given in Appendix A, and an implementation
can be found on the \href{https://www.github.com/xentlabs/s-xgl}{https://www.github.com/xentlabs/s-xgl}
repository.

\subsection{Examples of Xent Games}\label{subsec:examples-of-xent-games}

In this subsection, we review a number of Xent Games as a means of
illustrating the potential of the space. While many Xent Games can
be hand-designed for various purposes (see \cite{hongler-emil} for
a few), our thesis is that the means towards achieving general capabilities
is to rely on the automated design of such games. 

\subsubsection{Pre-Training Game}\label{subsec:pre-training-game}

The classical pre-training objective, i.e. the minimization of $\mathrm{xent}_{\mathcal{M}}$
is a canonical Xent Game. This game is in some sense the simplest:
there is no move elicited from $\mathcal{M}$, just a reward $-\mathrm{xent}_{\mathcal{M}}\left(x\right)$
where $x$ is a random string loaded from a data model. 

Variants of the basic pre-training objective, as e.g. the multi-token
prediction objective \cite{deepseek-v3} (e.g. asking a model to predict
the value for $10$ tokens ahead, without seeing the $9$ ones coming
before) can naturally be represented as Xent Games.

\subsubsection{Reinforcement Learning as Pre-Training (RLP)}\label{subsec:reinforcement-learning-as-pre-training}

In \cite{hatamizadeh-akter-prabhumoy-kautz-et-al}, a certain game,
called RLP, is considered (inspired by a similar game \cite{reinforcement-pretraining},
called RPT). The RLP task is in fact a Xent Game. Changing the notation
compared to the paper to fit the Xent Game description, and removing
a number of details that are not relevant to our discussion, the game
is the following: given a \emph{game map} $x_{<t},x_{t}$ taken from
a dataset (i.e. where $x_{t}$ follows the tokens of $x_{<t}$), the
game elicits a string $c$ from the player $\mathcal{M}$ to improve
the prediction of $x_{t}$ given $x_{<t}$ and $c$: the reward of
$\mathcal{M}$ is $-\mathrm{xent}_{\mathcal{M}}\left(x_{t}|x_{<t},c\right)$.
This is an interesting example of a game where the judge and the player
are the same model. 

\subsubsection{Distillation and Self-Distillation Games}

Distillation and self-distillation can both be viewed as instances
of Xent Games. 

In online distillation \cite{thinking-machines-on-policy,agarwal-on-policy-distillation},
we have a judge (teacher) model $\mathcal{J}$ that we want to distill
into a player (student) model $\mathcal{M}$. The associated Xent
Game uses the contrastive objective across models (Section \ref{subsec:Contrastive-objective}).
Concretely, we provide to $\mathcal{M}$ a context $x$ (possibly
void), and we elicit some answer $c$ from $\mathcal{M}$. The reward
for $\mathcal{M}$ is then $\mathrm{xent}_{\mathcal{M}}(c\mid x)-\mathrm{xent}_{\mathcal{J}}(c\mid x)$.
Maximizing this reward corresponds to minimizing a reverse-KL objective
between the (student) model $\mathcal{M}$ and the (teacher) model
$\mathcal{J}.$

In the articles \cite{hubotter-et-al26-self-distillation,shenfeld-self-distillation},
the authors do not consider the contrastive objective across models
but the counterfactual thinking one (Section \ref{subsec:Counterfactual-thinking}),
where the extra information comes either from an exact proof, or feedback
derived from previous attempts. A simplified Xent Game associated
with \cite{hubotter-et-al26-self-distillation} can be described as
follows. The main player $\mathcal{M}$ is cloned to obtain a judge
$\mathcal{J}$ used to provide reward. We provide $x$ to $\mathcal{M}$,
and elicit a continuation $c$ from $\mathcal{M}$. Next, conditioned
on $x,c$ we prompt $\mathcal{J}$ for feedback and obtain $f$. The
player then receives the reward $\mathrm{xent}_{\mathcal{M}}(c\mid x)-\mathrm{xent}_{\mathcal{J}}(c\mid x,f)$. 

\subsubsection{Prompt Games }\label{subsec:prompt-games}

Prompt games are our first concrete illustration of the idea of using
implicit knowledge to improve the explicit behavior. Inspired by Example
\ref{exa:inverse-prompting}, the basic reverse prompting game is
as follows. Consider a fixed judge $\mathcal{J}$; a \emph{game map}
consists of a text $s$. The player must find a prefix $t$ such that,
according to $\mathcal{J}$, $s$ becomes likely after $t$: we reward
the player with $-\mathrm{xent}_{\mathcal{J}}\left(s|t\right)$. This
game is related to a number of tasks, including creative summarization
or jailbreaking. Generalization follows from Example \ref{exa:common-explanations}
by taking as game map a set of texts $s_{1},\ldots,s_{n}$: the player
must find $t$ that makes them jointly more likely (for $\mathcal{J}$),
while remaining relatively unexpected from each of them individually.
Depending on the implementation details (constraints, prompting, regularization...),
variants of this inverse prompting game have different interpretations,
as discussed below. 
\begin{itemize}
\item Creative summarization. If we add constraints that prevent token copying
(length constraint, no common words, ... ), a regularization term
such as $-\mathrm{xent}_{\mathcal{J}}\left(t\right)/2$ (which is
equivalent to taking $-2\mathrm{xent}_{\mathcal{J}}\left(s|t\right)-\mathrm{xent_{\mathcal{J}}}\left(t\right)$)
to obtain a well-formed text $t$, and insert some sandwich prompt
between $s$ and $t$, then the optimal prefix $t$ can be considered
as a summary of $s$. 
\item Prompt injection. Prompt injection tasks can be formulated as Xent
Games. A simple case is the static objective: find a prompt (or injected
data) that causes a model $\mathcal{J}$ to generate a fixed output
$s$, independent of the user's instructions $i$ and data $d$. To
evaluate the gap between the response generated by the LLM $\mathcal{J}$
and the target answer $s$, the authors of \cite{automatic-universal-prompt-injection-attacks}
propose using a xent loss, e.g. $\mathrm{xent}_{\mathcal{J}}(s\mid i,d,x)$.
In the corresponding Xent Game, another model (the attacker) proposes
$x$ (conditioning on $(i,d)$ or not), to maximize the expected reward
over a training set of instructions-data pairs. 
\item Defensive variants. Defensive methods such as DataSentinel \cite{data-sentinel}
can be viewed as a two-player Xent Game, whose goal is to prevent
injected tasks. Two models, an attacker $\mathcal{A}$ and a detector
$\mathcal{D}$, compete; this leads to a game-theoretic min-max optimization
problem. The attacker outputs contaminated prompts or data (given
some instruction and data to process) so that a judge LLM $\mathcal{J}$
executes an injected task instead of the intended task, while also
avoiding to be detected. The detector receives a fixed detection instruction
$s_{d}$ and the (possibly contaminated) prompt/data and outputs a
secret detection string which serves as marker. This marker should
be present when the input is clean, and absent from the detector's
output if contaminated. In this setup, all losses are cross-entropy
terms. For example, detection is quantified by the cross entropy $\mathrm{xent}_{\mathcal{G}}(k\mid s_{d},x)$
of the marker $k$ given the detection instruction $s_{d}$ concatenated
with the potentially contaminated data $x$. 
\end{itemize}

\subsubsection{Approximate Verifiable Reward Games}\label{subsec:approximate-verifiable-reward-games}

When verification and scoring is easier than construction, a judge
model can be used to enforce legality of player moves and score plausibility.
As explained in Section \ref{subsec:verification-vs-construction},
such a well-trained judge, even if imperfect, can provide a continuous
signal: reward in a verifiable task is thus turned into continuous
feedback that can drive parameter optimization. 

One may worry that the player will discover ways to exploit the judge
biases rather than solve the intended verifiable task. While it is
always possible to add an external reward for truly solving the problem,
the Xent Games usage relies on transfer across games: if reward hacking
is too easy, we conjecture that the model will not learn new generalizable
capabilities, and such a game will be filtered out when games are
selected in the curriculum. 

Examples of such verifiable reward games include chess and mathematical
proofs (as discussed in \cite{hongler-emil}):
\begin{itemize}
\item Chess. Games like chess can be used to define Xent Games, based on
a judge model which is strong enough to 1. recognize legal moves and
2. evaluate positions. A game map in this case is a chess position.
Two players alternately output moves, ensure functions enforce each
move's validity. After a fixed number of steps, or at termination,
the score is defined from the judge's cross-entropy difference between
statements such as ``white is winning'' and ``black is winning''. 
\item Mathematical proofs. Let us again assume that a judge model $\mathcal{J}$
is strong enough to 1. estimate the correctness of short formal proofs
(for example, in Lean) and 2. evaluate the plausibility of mathematical
proof sketches written in English. We can adapt the debate protocol
of \cite{cghcl-i,cghcl-ii}, where two agents explore a tree of proof
sketches: the prover wants to maximize the plausibility of the proof
by adding details, while the skeptic wants to minimize it by asking
for details. The economic rewards of the SPRIG protocol are then replaced
by xent terms using the judge model, defining a zero-sum Xent Game.
\end{itemize}

\subsection{Xent Game Space Properties}\label{subsec:xent-game-space-properties}

As argued in Section \ref{subsec:examples-of-xent-games}, the space
of Xent Games $\mathcal{G}$ forms a vast family, eliciting numerous
skills highlighted in other papers. In this subsection, we argue that
this space possesses a number of good properties that make it suitable
for automated exploration towards building generally capable agents. 

\subsubsection{Closure under Axioms}

The first natural property of the Xent Game space is its closure under
natural axioms. In \cite{hongler-emil}, the O-XG space can be shown
to naturally emerge from a small collection of game-theoretic and
compositional axioms. The S-XG space can be characterized very similarly:
\begin{claim}
\label{claim:sxg-characterization}If a space of text games on token
string space (with the operations described in Section \ref{subsec:token-string-space})
contains the reverse prompt game (Section \ref{subsec:prompt-games})
and is stable under compositionality (we can append the instructions
of one game to another), zero-summing (a player's loss can be a player's
gain and vice versa) and adversarial rescalings (we can turn \emph{reward}
statements into \emph{ensure} ones), then it must contain the space
of S-XG.
\end{claim}

\begin{rem}
In \cite{hongler-emil}, adversarial rescaling is formulated in terms
of unbounded multipliers, leading to hard constraints; for the sake
of training models (rather than evaluating them), it is better to
work with soft constraints (with bounded multipliers); the principle
is exactly the same. Similarly, the set of allowed operations on the
string space in S-XG is slightly different from that of O-XG, but
the principle is exactly the same, and if we use the operations described
in Section \ref{subsec:token-string-space}, we obtain the Xent Games
described in Section \ref{subsec:sxg-def-and-basic-runtime-instruction}. 
\end{rem}

\subsubsection{Web-of-Games Assumption and Curriculum Universality}\label{subsec:web-of-games}

A second important conjectural property of the Xent Game space outlined
in \cite{hongler-emil} concerns the transfer value: informally, the
idea is that it is relatively easy, given a set of games, to ``discover''
new, related games, allowing one to navigate in the space of games,
and to grow the skills of a model. In the context of curriculum learning,
this motivates the following approximation claim about curricula made
of general verifiable reward environments (i.e that are not necessarily
Xent Games): 
\begin{claim}
\label{claim:web-of-games}Let $E$ be a skill evaluation. Suppose
there exists a curriculum of verifiable reward environments/tasks
$\mathbb{G}_{1},\ldots,\mathbb{G}_{n}$ such that training on this
curriculum brings any model $\mathcal{M}$ in a family $\mathfrak{M}$
to a target skill level $\Lambda_{\mathcal{M}}$ on $E$. Then, for
any $\epsilon>0$, there exists
\begin{itemize}
\item a judge model $\mathcal{M}_{\mathcal{J}}$,
\item a curriculum of Xent Games $G_{1},\ldots,G_{m}$ of Xent Games that
is algorithmically computable from $\mathbb{G}_{1},\ldots,\mathbb{G}_{n}$, 
\end{itemize}
such that training any $\mathcal{M}\in\mathfrak{M}$ on $G_{1},\ldots,G_{m}$
achieves skill level at least $\Lambda_{\mathcal{M}}-\epsilon$ on
$E.$ Moreover, suppose that the curriculum $\mathbb{G}_{1},\cdots,\mathbb{G}_{n}$
is computed by a greedy optimization process, then the curriculum
of Xent Games $G_{1},\ldots,G_{m}$ can also be computed using a greedy
optimization process. 
\end{claim}

Taken at face value, this claim is not very surprising: if (for instance),
we have an environment with binary outputs we can demand that $\mathcal{M}_{\mathcal{J}}$
be strong enough to evaluate the games of $\mathbb{G}_{1},\ldots,\mathbb{G}_{n}$
and approximate (naively) the scores by looking at the xent of ``the
output of $\mathbb{G}_{j}$ is $0$'' versus the xent of ``the output
of $\mathbb{G}_{j}$ is $1$''. This is similar to the ideas of \ref{subsec:approximate-verifiable-reward-games}. 

Similarly to the universality theorems for neural networks, the means
to justify Claim \ref{claim:web-of-games} are not necessarily directly
informative of practical/relevant uses of Xent Games; rather, the
statement's point is that ``nothing is missing'' from the space.

In Section \ref{sec:cognitive-training}, we leverage these assumptions
to formulate cognitive training as a meta-sampling process over the
space of Xent Games. 

\subsection{Game Training and Transfer}\label{subsec:game-training-and-transfer}

The purpose of Xent Games is to provide a suitable environment for
robust learning. In this subsection, we introduce the key notion of
transfer value between games, which quantifies how training under
a scheme $\Phi$ on one game affects the performance on another. 

\subsubsection{Training Scheme}\label{subsec:training-scheme}

We consider a fixed training scheme $\Phi$, which defines, for any
Xent Game $G$, a map $\mathcal{M}\mapsto\Phi_{G}\mathcal{M}$ which
returns the model obtained by training $\mathcal{M}$ on $G$ for
\emph{one run} of the game. A game with multiple training steps should
be thought of in our framework as being made of many concatenated
copies of a game $G\oplus G\oplus\cdots\oplus G$.
\begin{rem}
For the sake of training, it can be considered that training steps
are performed when a clearing ``end-of-game'' instruction is called. 
\end{rem}

We assume that the training scheme $\Phi$ is invariant by rescaling
(if $\alpha G$ denotes the games obtained by multiplying all scores
by $\alpha$, then training is identical and we obtain $\Phi_{\alpha G}=\Phi_{G}$;
in other words, the training is invariant by a rescaling of units
of score). This is the case for GRPO-style training based on centered
and normalized rewards, as well as for algorithms such as \cite{hongler-emil-renard-gabriel,hongler-asani-gabriel-emil-renard}. 
\begin{rem}
\label{rem:important-theoretical-possibility} Importantly: the space
of Xent Games is not suitable to perform arbitrary rescalings of game
scores by $\alpha>0$. We \emph{could theoretically allow} the S-XGL
language to support a global rescaling of scores for each game, but
this would not change anything about what can be achieved in terms
of model training (as a result of the scale invariance of the training).
However this theoretical possibility (which reflects a training algorithm)
is important for our derivation (see Section \ref{subsec:meta-objective-derivation}) 
\end{rem}

\subsubsection{Transfer Value}\label{subsec:transfer-value}

A central quantity in our approach is the transfer value between games.
Informally, the transfer value $\mathcal{T}_{G}^{\mathcal{M}}\left(H\right)$
encodes ``how much in expectation does training on $G$ teaches $\mathcal{M}$
about how to play $H$''. We denote the expected score of $\mathcal{M}$
on $H$ by 
\[
\mathrm{S}_{\mathcal{M}}\left(H\right):=\mathbb{E}\left[\mathrm{Score}_{H}\left(\mathcal{M}\right)\right].
\]
 The transfer value is defined as follows. 
\begin{defn}
\label{def:transfer}For two games $G$ and $H$, the transfer value
$\mathcal{T}_{G}^{\mathcal{M}}\left(H\right)$ from $G$ to $H$ for
$\mathcal{M}$ is 
\[
\mathcal{T}_{G}^{\mathcal{M}}\left(H\right):=\mathrm{S}_{\Phi_{G}\left[\mathcal{M}\right]}\left(H\right)-\mathrm{S}_{\mathcal{M}}\left(H\right).
\]
More generally, if $E$ is an external evaluation and $G$ is a Xent
Game, we define 
\[
\mathcal{T}_{G}^{\mathcal{M}}\left(E\right):=\mathrm{S}_{\Phi_{G}\left[\mathcal{M}\right]}\left(E\right)-\mathrm{S}_{\mathcal{M}}\left(E\right),
\]
where $\mathrm{S}_{\cdot}(E)$ is the expected reward on $E$. 
\end{defn}

\begin{rem}
While these notions are well-posed theoretically, estimating them
in practice may require a large number of samples or other techniques
(such as the upcoming \cite{hongler-emil-renard-gabriel,hongler-asani-gabriel-emil-renard}). 
\end{rem}

A useful scaling result is the following:
\begin{rem}
\label{rem:game-scaling}For games $G,H$ we have the scaling relations: 
\begin{itemize}
\item $\mathcal{T}_{G}^{\mathcal{M}}\left(H\oplus H\right)=2\mathcal{T}_{G}^{\mathcal{M}}\left(H\right)$:
this is as we are purely evaluating on $H$ (twice) and the model
is not learning between steps.
\item For a small game $G$, we have $\mathcal{T}_{G\oplus G}^{\mathcal{M}}\left(H\right)\approx2\mathcal{T}_{G}^{\mathcal{M}}\left(H\right)$;
training twice on it yields (in expectation) approximately twice the
increase of performance on $H$ (this is only an approximation, as
if we keep repeating a game enough, we may see diminishing returns).
\end{itemize}
\end{rem}

In Section \ref{sec:cognitive-training} below, transfer value will
be key to:
\begin{itemize}
\item estimating the internal relevance of candidate new games;
\item estimating the novelty brought by new games compared to previously
selected games. 
\end{itemize}

\subsubsection{Positive Correlation}\label{subsec:positive-correlation}

A fundamental assumption upon which cognitive learning rests is that
we can work with games that allow us to grow a curriculum constructively
in a monotone way, i.e. without ``needing to make a step back for
every two steps forward''. 

For two Xent Games $G_{1,}G_{2}\in\mathcal{G}$, we say that they
are positively correlated (relative to a model $\mathcal{M}$) if
$\mathcal{T}_{G_{i}}^{\mathcal{M}}\left(G_{j}\right)>0$ for $i,j\in\left\{ 1,2\right\} $:
informally speaking, this simply means that learning one doesn't make
$\mathcal{M}$ worse at the other. To grow a curriculum of games in
$\mathcal{G}$, we want to be able to pick at any time (and then optimize
with respect to the meta-objective, see Section \ref{sec:cognitive-training})
new games that are positively correlated with respect to all the past
games of the curriculum, as measured at appropriate times (see also
Section \ref{subsec:meta-objective-derivation} below for a justification
of the times): if we have built a curriculum $G_{0},\ldots,G_{k-1}$
yielding the models $\mathcal{M}_{1},\ldots,\mathcal{M}_{k}$, we
want to be able to pick a next game $G_{k}$ such that (\emph{a minima})
for each $j<k$, we have:
\begin{align}
\mathcal{T}_{G_{j}}^{\mathcal{M}_{j}}\left(G_{k}\right) & >0,\label{eq:old-to-new-nonneg}\\
\mathcal{T}_{G_{k}}^{\mathcal{M}_{k}}\left(G_{j}\right) & >0.\label{eq:new-to-old-nonneg}
\end{align}

\begin{defn}
We define $\mathcal{G}_{k}^{+}\subset\mathcal{G}$ to be the set of
Xent Games that are positively correlated to $G_{<k}$.
\end{defn}

\begin{rem}
An undesirable outcome that we want to avoid is that we end up in
a \emph{transfer cul-de-sac}, i.e that we are not able to further
pick any new game (even a copy of an old one) satisfying (\ref{eq:old-to-new-nonneg}
and \ref{eq:new-to-old-nonneg}). In case an otherwise useful game
happens to negatively correlate with one (or a few) old games, a solution
can be simply to append to it a number of copies of these old games
to prevent regression on them (and hence negative correlation). 
\end{rem}

\subsubsection{Relevant Skill Discovery}

An important desirable feature of a curriculum of games $\left(\mathcal{G}_{k}\right)_{k\geq0}$,
which motivates cognitive training is that of \emph{relevant skill
discovery}:
\begin{defn}
\label{def:relevant-skill-discovery}We say that a sequence of games
$\left(G_{k}\right)_{k\geq0}$ achieves \emph{relevant skill discovery}
on $\mathcal{M}$ if, when training a sequence of models $\left(\mathcal{M}_{k}\right)_{k}$,
with $\mathcal{M}_{0}=\mathcal{M}$ and $\mathcal{M}_{k+1}$ obtained
from $\mathcal{M}_{k}$ by training on $G_{k},$we achieve the following:
\begin{itemize}
\item The process keeps discovering new skills, i.e. finding new games that
are not completely covered (in the sense of transfer value) by any
$G_{<k}$ for a fixed $k$.
\item It keeps growing on already discovered skills, i.e. for any $G_{j}$
discovered, performance on $G_{j}$ keeps improving as $k$ increases. 
\end{itemize}
\end{defn}

\section{Cognitive Training}\label{sec:cognitive-training}

In the previous section, we described the Xent Game Space $\mathcal{G}$,
which we postulate provides a good coverage of the learnable skills
of an LLM (see Section \ref{subsec:web-of-games}). Our goal is now
to provide an algorithm to construct a curriculum of games in $\mathcal{G}$
for cognitive training, i.e. to build a ``flywheel'' that continuously
(and autonomously) provides a stream of games that improve the model's
general skills. 
\begin{rem}
\label{rem:single-model-trained}In this section, we focus on the
training of \emph{a single model} $\mathcal{M}$, although we may
use other frozen (non-trained) models as judges or data streamers.
Simultaneously training a family of models can be an interesting way
to make cognitive training more efficient, but goes beyond the scope
of this paper. 
\end{rem}

In this section, we describe a meta-objective to gauge the quality
of a stream of games, based on elaborations of the evolution-based
ideas of \cite{hongler-emil}. Given an initial model $\mathcal{M}=\mathcal{M}_{0}$,
the goal is to evolve it into a sequence of models $\mathcal{M}_{1},\mathcal{M}_{2},\ldots$
that are more and more generally capable, where each update $\mathcal{M}_{k}\mapsto\mathcal{M}_{k+1}$
is obtained by training on a game $G_{k}$. 

\subsection{Cognitive Training as a Greedy Algorithm}

The central question is \emph{how to choose the next game in the cognitive
training curriculum}. We propose to approach this as a greedy optimization
problem, focusing on picking the \emph{next game} to grow the curriculum.
\begin{rem}
\label{rem:greedy-algorithm}We focus on picking a greedy algorithm,
as the task of optimizing a true value function on all curricula is
likely completely intractable as the curriculum grows. We postulate
that this greedy process is sufficient towards discovering new skills
and growing towards general capabilities. 
\end{rem}

At each step, the goal is to find the best (most valuable) game that
allows us to evolve the current model towards higher ``cognitive''
abilities, i.e. that has acquired ``the most valuable new skills''
between steps $k$ and $k+1$. 

The \emph{meta-game at stage $k$} consists in tasking a meta-sampler
model $\mathcal{M}_{\mathcal{S}}$ with outputting a Xent Game $G_{k}$
(in S-XGL): the reward of $\mathcal{M}_{\mathcal{S}}$ is evaluated
by a meta-objective $G_{k}\mapsto\mathcal{O}\left(G_{k}\right)$,
which depends on $G_{<k}$ and $\mathcal{M}_{\leq k}$. 

The question of determining a principled form for $\mathcal{O}$ is
a priori a very underspecified question: assuming that the goal of
the meta-game is to ascribe a value to games, this raises the question
of whether there is a meta-meta-objective, and so on and so forth...
it is not clear that such a problem can be resolved in a constructive
fashion, i.e. without relying upon an infinite number of assumptions
or parameters. \footnote{This situation echoes a bit a similar situation in field theory, where
a Lagrangian may define an action of govern fields in space-time via
e.g. a least-action principles, while there is no least-action principle
that a Lagrangian would itself satisfy a priori. }

\subsection{Meta-Objective Formula}\label{subsec:meta-objective-formula}

In Section \ref{subsec:meta-objective-derivation}, we show that (perhaps
surprisingly) an explicit determination of $\mathcal{O}$ from purely
qualitative principles (such as consistency principles) can be obtained,
substantially constraining the space of ``sound'' meta-objectives
to an explicit natural choice formula\footnote{A somewhat analogous situation appears in the so-called bootstrap
program in conformal field theory.}.

Given a history $G_{<k}$ and resulting models $\mathcal{M}_{\leq k}$,
and restricting the domain to the positively-correlated games $\mathcal{G}_{k}^{+}$
(see Section \ref{subsec:positive-correlation}), we find 
\begin{equation}
\mathcal{O}=\frac{qd+bp}{l},\label{eq:meta-objective}
\end{equation}
where we have
\begin{itemize}
\item the \emph{quality} term $q$ is given by $q\left(H\right)=\left(\sum_{j<k}\mathcal{T}_{H}\left(G_{j}\right)\right)^{1-\delta}$
for a diversity hyper-parameter $\delta\in\left[0,1\right]$
\item the \emph{diversity} term $d$ is given by $d\left(H\right)=\left(\mathcal{T}_{H}\left(H\right)/\sum_{j<k}\mathcal{T}_{G_{j}}\left(H\right)\right)^{\delta}$, 
\item the \emph{benchmark} term $b$ is given by $b\left(H\right)=\mathcal{T}_{H}\left(E\right)$
for an external benchmark metric term $E$,
\item the \emph{pressure} term $p\geq0$ is a hyper-parameter modulating
the importance of $qd$ versus $b$.
\item the \emph{length} term $l$ is the raw code length of a game written
in S-XGL. 
\end{itemize}
Note that despite a linearly growing number of terms in $k$ for $q$
and $d$ (suggesting quadratic scaling in $k$ to run Cognitive Training),
the meta-objective computation is only growing very slowly in $k$
(see Section \ref{subsec:scalability-of-meta-objective-computation}). 

The derivation is presented in Section \ref{subsec:meta-objective-derivation}
below; in particular, the key result is the derivation of the Internal
Meta-Objective Theorem derived in Section \ref{subsec:internal-meta-objective-derivation}. 
\begin{rem}
The hyper-parameters $\delta$ and $p$ can in principle depend on
the training step; if they do, their formulation is however constrained
by similar principles to those used in the derivation below. 
\end{rem}

\subsection{Meta-Objective Principles}\label{subsec:meta-objective-derivation}

In this section, we derive an expression for the meta-objective $\mathcal{O}$
aimed at assessing the value of a new game $G_{k}$, given a past
curriculum $G_{<k}$ of ``old games'' used to train the model sequences
$\mathcal{M}_{<k}$. 

\subsubsection{General Principles}\label{subsec:general-principles}

The following line of reasoning leads us to a principled expression
for $\mathcal{O}$ 
\begin{itemize}
\item The goal of $\mathcal{O}$ is to measure the value that a new game
brings towards training the next model.
\item The $\mathcal{O}$ value combines \emph{internal} and \emph{external}
relevance terms, denoted by $\mathcal{I}$ and $\mathcal{E}$ respectively.
\item Both $\mathcal{I}$ and $\mathcal{E}$ terms measure the ``useful
work performed'' if we train the model $\mathcal{M}_{k}$ on $H$. 
\item The $\mathcal{I}$ and $\mathcal{E}$ terms of a game $H$ are to
be normalized by its code length $l$ in number of tokens. 
\item All terms should be evaluated in an \emph{online} fashion, i.e. not
allowing the re-training of archives $\mathcal{M}_{<k}$.
\item The $\mathcal{E}$ term depends on the difference in an external benchmark
metric $E$ if we train $\mathcal{M}_{k}$ on $H$.
\end{itemize}
We thus obtain the following expression for $\mathcal{O}$ of the
form, 
\[
\mathcal{O}=\frac{\mathcal{I}+\mathcal{E}}{l},
\]
in what follows, we will provide a principled derivation for $\mathcal{I}$
and $\mathcal{E}$.

\subsubsection{Internal Meta-Objective Derivation}\label{subsec:internal-meta-objective-derivation}

To provide a principled derivation for the internal relevance $\mathcal{I}$
of the meta-objective $\mathcal{O}$ is quite nontrivial: the problem
is a priori heavily underspecified. Assume again that we have obtained
$\mathcal{M}$ by training on $G_{<k}$ leading to the sequence $\mathcal{M}_{\leq k}$. 
\begin{enumerate}
\item The term $\mathcal{I}$ depends on the past games $G_{<k}$ via the
following transfer values:
\begin{enumerate}
\item The ``new-to-old'' transfer values $\mathcal{T}_{H}^{\mathcal{M}_{k}}\left(G_{j}\right)$
for $j<k$: this is the core of the \emph{quality} measure.
\item The ``old-to-new'' transfer values $\mathcal{T}_{G_{j}}^{\mathcal{M}_{j}}\left(H\right)$
for $j<k$: this is the core of the \emph{diversity} measure.
\item The ``new-to-new'' transfer value $\mathcal{T}_{H}^{\mathcal{M}_{k}}\left(H\right)$:
this is a term that is useful for normalization. 
\end{enumerate}
\item The term $\mathcal{I}$ should be factored as 
\[
\mathcal{I}=qd,
\]
where the \emph{quality} $q$ depends on $\left(\mathcal{T}_{H}^{\mathcal{M}_{k}}(G_{j})\right)_{j<k}$
and the \emph{diversity} $d$ depends on $\left(\mathcal{T}_{G_{j}}^{\mathcal{M}_{j}}(H)\right)_{j<k}$
and on $\mathcal{T}_{H}^{\mathcal{M}_{k}}\left(H\right)$; using a
product (rather than a sum) ensures the optimization of \emph{both}
quantities. 
\item The quality and diversity terms are invariant under \emph{history
fusion}: if for $0<j<k$, we replace the steps 
\[
\mathcal{M}_{j-1}\underset{G_{j-1}}{\longrightarrow}\mathcal{M}_{j}\underset{G_{j}}{\longrightarrow}\mathcal{M}_{j+1}
\]
by the training step (where $\mathrm{\left[none\right]}$ is an ``idle''
game, performing nothing)
\[
\mathcal{M}_{j-1}\underset{G_{j-1}\oplus G_{j}}{\longrightarrow}\mathcal{M}_{j+1}\underset{\left[\mathrm{none}\right]}{\longrightarrow}\mathcal{M}_{j+1},
\]
then values of $q$ and $d$ for any $G_{k}$ are left unchanged. 
\item From this, we deduce that $q$ must only depend on the sum of transfer
values, i.e. must be of the form
\begin{align*}
q\left(H\right)= & f_{q}\left(\sum_{j<k}\mathcal{T}_{H}^{\mathcal{M}_{k}}\left(G_{j}\right)\right),\\
d\left(H\right)= & f_{d}\left(\sum_{j<k}\mathcal{T}_{G_{j}}^{\mathcal{M}_{j}}\left(H\right),\mathcal{T}_{H}^{\mathcal{M}_{k}}\left(H\right)\right).
\end{align*}
 
\item Using the ``theoretical rescaling invariance idea'' (see Remark
\ref{rem:important-theoretical-possibility}), we get the following:
we must have rescaling invariance for $q\left(H\right)$, i.e. $q\left(\alpha H\right)=q\left(H\right)$
for $\alpha>0$. From this idea, we also get that $\mathcal{I}\left(\alpha H\right)=\mathcal{I}\left(H\right)$
(since the value of a game can only depend on its effect for training
models and training is invariant under rescaling). From there, we
find that $d\left(\alpha H\right)=d\left(H\right)$ as well. Since
$\mathcal{T}_{G_{j}}^{\mathcal{M}_{j}}\left(\alpha H\right)=\alpha H$
and $\mathcal{T}_{\alpha H}^{\mathcal{M}_{k}}\left(\alpha H\right)=\alpha\mathcal{T}_{H}^{\mathcal{M}_{k}}\left(H\right)$,
we find that $f_{d}$ must be a function of the \emph{ratio} of its
two arguments only, i.e. be of the form
\[
d\left(H\right)=F_{d}\left(\frac{\mathcal{T}_{H}^{\mathcal{M}_{k}}\left(H\right)}{\sum_{j<k}\mathcal{T}_{G_{j}}^{\mathcal{M}_{j}}\left(H\right)}\right)
\]
 for an increasing bijective function $F_{d}:\mathbb{R}_{+}\to\mathbb{R}_{+}$. 
\item From Remark \ref{rem:game-scaling}, we have that $\mathcal{T}_{H\oplus H}^{\mathcal{M}_{k}}\left(H\oplus H\right)\approx4\mathcal{T}_{H}^{\mathcal{M}_{k}}(H)$,
and requiring $\mathcal{I}\left(H\oplus H\right)\approx2\mathcal{I}\left(H\right)$
(the usefulness of training on a game twice is approximately twice
the usefulness of the game), we obtain $qd\left(H\oplus H\right)\approx2qd\left(H\right)$
and assuming that $f_{q}$ and $F_{d}$ must be homogeneous increasing
functions, we find that their powers sum up to $1$, yielding 
\[
\left(\sum_{j<k}\mathcal{T}_{H}^{\mathcal{M}_{k}}\left(G_{j}\right)\right)^{1-\delta}\left(\frac{\mathcal{T}_{H}^{\mathcal{M}_{k}}\left(H\right)}{\sum_{j<k}\mathcal{T}_{G_{j}}^{\mathcal{M}_{j}}\left(H\right)}\right)^{\delta},
\]
for some \emph{diversity parameter} $\delta\in\left[0,1\right]$. 
\end{enumerate}
From the above, we obtain the following: 
\begin{thm*}[Internal Meta-Objective Theorem]
A consistent internal meta-objective (according to the principles
outlined above) $\mathcal{I}$ must be of the form 
\[
\mathcal{I}=qd,
\]
where 
\begin{align*}
q\left(H\right) & =\left(\sum_{j<k}\mathcal{T}_{H}^{\mathcal{M}_{k}}\left(G_{j}\right)\right)^{1-\delta},\\
d\left(H\right) & =\left(\frac{\mathcal{T}_{H}\left(H\right)}{\sum_{j<k}\mathcal{T}_{G_{j}}^{\mathcal{M}_{j}}\left(H\right)}\right)^{\delta},
\end{align*}
where $\delta\in\left[0,1\right]$ is a diversity hyper-parameter.
\end{thm*}

\subsubsection{External Meta-Objective Derivation}\label{subsec:external-meta-objective-derivation}

We can apply a similar idea to form a natural external meta-objective: 
\begin{enumerate}
\item Similarly to $\mathcal{I}$, the term $\mathcal{E}$ should satisfy
$\mathcal{E}\left(H\oplus H\right)\approx2\mathcal{E}\left(H\right)$.
\item As a result, it is natural to postulate that 
\[
\mathcal{E}=bp,
\]
where $b=\mathcal{T}_{H}^{\mathcal{M}_{k}}\left(E\right)$, with $E$
being the given external benchmark metric, where $p$ is a hyper-parameter.
\end{enumerate}

\subsubsection{On the Hyper-Parameters $\delta$ and $p$}\label{subsec:on-hyper-parameters}

At every step, the two hyper-parameters $\delta>0$ and $p\geq0$
are in principle allowed to vary; a principled derivation of formulae
for those (i.e. how they should vary in time) is a priori non-trivial
and lies beyond the scope of the present work.

A simple choice for both is to take them as constants or to rely on
some scaling with respect to the time $T=\sum_{j<k}l\left(G_{j}\right)$
(note that this expression for $T$ is constrained by history fusion). 

\subsubsection{Remarks on the Derivation}\label{subsec:remarks-on-derivation}

In Sections \ref{subsec:general-principles}--\ref{subsec:external-meta-objective-derivation},
a principled derivation for $\mathcal{O}$ is provided. A number of
remarks are in order.
\begin{itemize}
\item The naming for the \emph{quality} and \emph{diversity }terms $q$
and $d$ comes from the Quality-Diversity line of research in Artificial
Life (see e.g. \cite{quality-diversity,etcheverry-chan-moulin-frier-oudeyer}),
but it is interesting to see that the derivation is made out of first
principles only (as opposed to a modeling attempt).
\item The derivation can be viewed as the result of the answer to the question:
what are the most natural objects we are able to use, and what is
a consistent way to combine them? It is notable that a few natural
constraints can narrow down the space of consistent meta-objectives
to a few hyper-parameters. 
\item The derivation is not restricted to Xent Games: this is discussed
in Section \ref{subsec:taking-a-step-back}.
\end{itemize}

\subsubsection{Scalability of Meta-Objective Computation}\label{subsec:scalability-of-meta-objective-computation}

While we assume access to the whole set of archives $\mathcal{M}_{j}$
for $j\leq k$, it is notable that to compute the diversity denominator
sum $\sum_{j<k}\mathcal{T}_{G_{j}}^{\mathcal{M}_{j}}\left(H\right)$
only requires us to measure the difference of performance of $H$
on $\mathcal{M}_{k}$ and $\mathcal{M}$: summing the Definition \ref{def:transfer},
we get a telescoping sum, where the internal terms cancel. 

Thus to perform Cognitive Training (Section \ref{subsec:cognitive-training-algorithm}),
we only need to store the \emph{latest archive} of $\mathcal{M}_{k}$.
Similarly, to estimate the other sum $\sum_{j<k}\mathcal{T}_{H}^{\mathcal{M}_{k}}\left(G_{j}\right)$
only obviously requires access to $\mathcal{M}_{k}$, but also can
be computed exactly in a time \emph{linear in the number} \emph{of
different games} used (and furthermore is easily parallelizable). 

\subsection{Cognitive Training Algorithm}\label{subsec:cognitive-training-algorithm}

In Section \ref{subsec:meta-objective-derivation}, we have derived
the explicit form for the meta-objective $\mathcal{O}$ that, given
an initial segment of the curriculum $G_{<k}$, ascribes value to
a new game $H\in\mathcal{G}_{k}^{+}$ (the space of Xent Games positively
correlated to $G_{<k}$, see Section \ref{subsec:positive-correlation}).
From this we first derive a cognitive training step:
\begin{defn}[Cognitive Training Step]
\label{def:cognitive-training-step}Given a step $k$ and a maximal
length $L$, and a sequence of games $G_{<k}$ used to train a sequence
of models $\left(\mathcal{M}_{k}\right)_{k\geq0}$ with $\mathcal{M}_{0}=\mathcal{M}$
and $\mathcal{M}_{k+1}$ being obtained by training $\mathcal{M}_{k}$,
optimize the meta-objective value $\mathcal{O}\left(H\right)$ on
games of length $l\leq L$. 
\end{defn}

\begin{rem}
Playing this game is a priori hard, as it involves an optimization
on the space of all games. As a result, playing can only be imperfect,
and can be hard to learn (see Section \ref{subsec:playing-the-meta-game}).
\end{rem}

From there, the cognitive training is defined as follows: 
\begin{defn}
\label{def:cognitive-training-algorithm}Given a meta-sampler $\mathcal{M}_{\mathcal{S}}$
(possibly with some prompt), the Cognitive Training consists in the
iteration of Cognitive Training Steps (see Definition \ref{subsec:cognitive-training-algorithm})
played by $\mathcal{M}_{\mathcal{S}}$.
\end{defn}

\begin{rem}
Note that a priori, the cognitive training algorithm is not a meta-game
itself (it has no single objective), being rather the result of iterated
play; it is rather a greedy optimization algorithm. 
\end{rem}

A toy model implementation of Cognitive Training, highlighting the
role of the diversity parameter $\delta$, can be found in \href{https://www.github.com/xentlabs/cognitive-training}{https://www.github.com/xentlabs/cognitive-training}.

\subsection{Toy Model for Cognitive Training}\label{subsec:toy-model-for-cognitive-training}

To give an intuition for the Cognitive Training process, we simulate
the greedy optimization of the meta-objective $\mathcal{O}$ on a
toy example:
\begin{itemize}
\item We assume the existence of a finite list of independent latent skills. 
\item Each game $G$ is associated with a skill vector that describes how
much each skill is required for playing $G$ (and, accordingly, how
much playing $G$ improves each skill). 
\end{itemize}
Modeling details are given in Appendix B. Figure \ref{fig:simulation}
shows the results of running the algorithm for 2048 steps, comparing
it with a baseline that randomly selects a game to train on at each
step. 
\begin{itemize}
\item The left panel shows the dependency of the training on $\delta$:
pure quality maximization ($\delta=0$) performs slightly worse than
random choice, whereas any sufficiently large $\delta$ leads to a
substantial improvement over the baseline. 
\item The middle panel additionally indicates faster learning with the greedy
algorithm throughout the whole training process. 
\item The right panel shows that the greedy algorithm uses only a small
subset of all 256 available games. 
\end{itemize}
\begin{figure}
\includegraphics[scale=0.25]{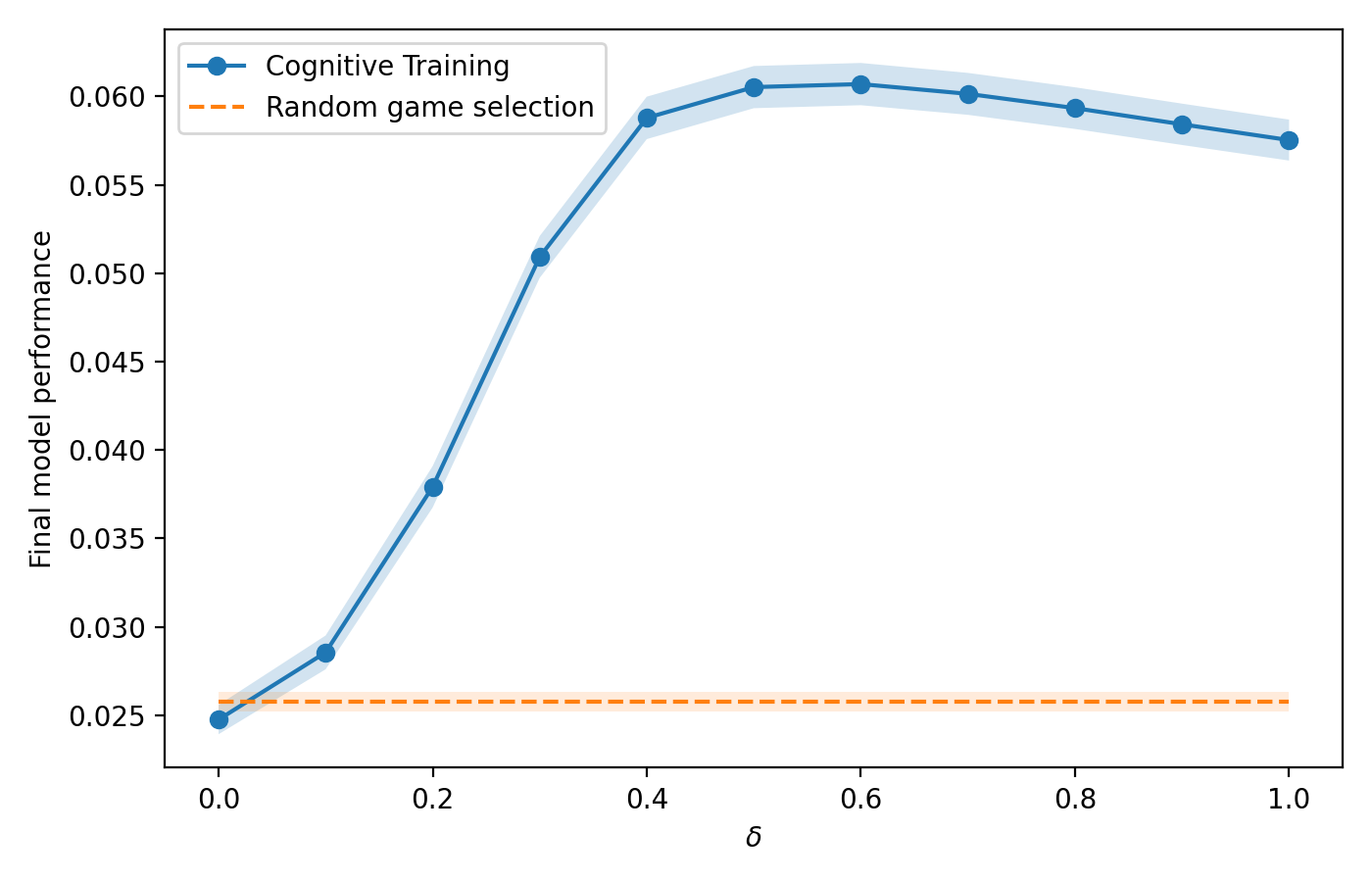} \includegraphics[scale=0.25]{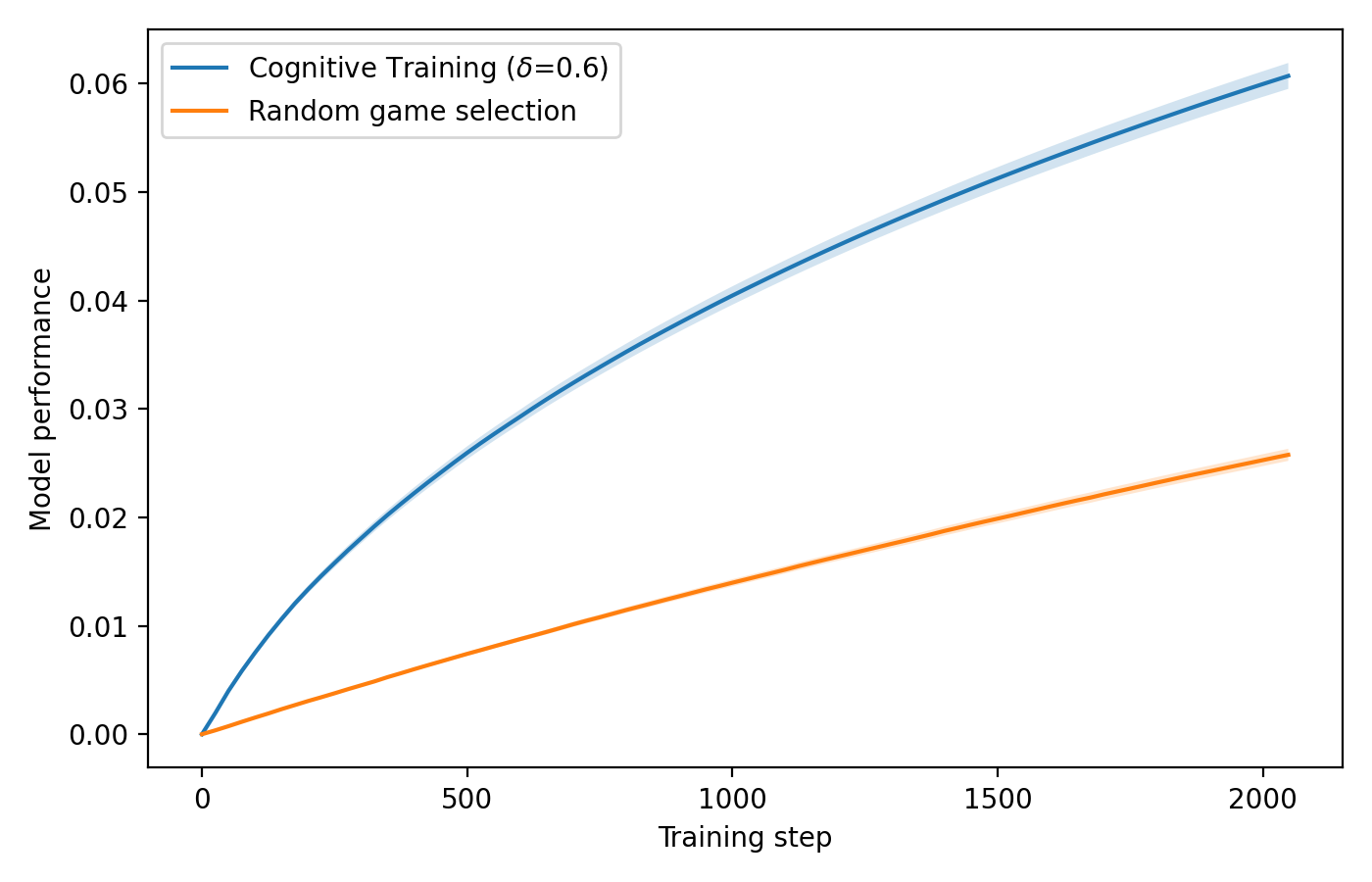}
\includegraphics[scale=0.25]{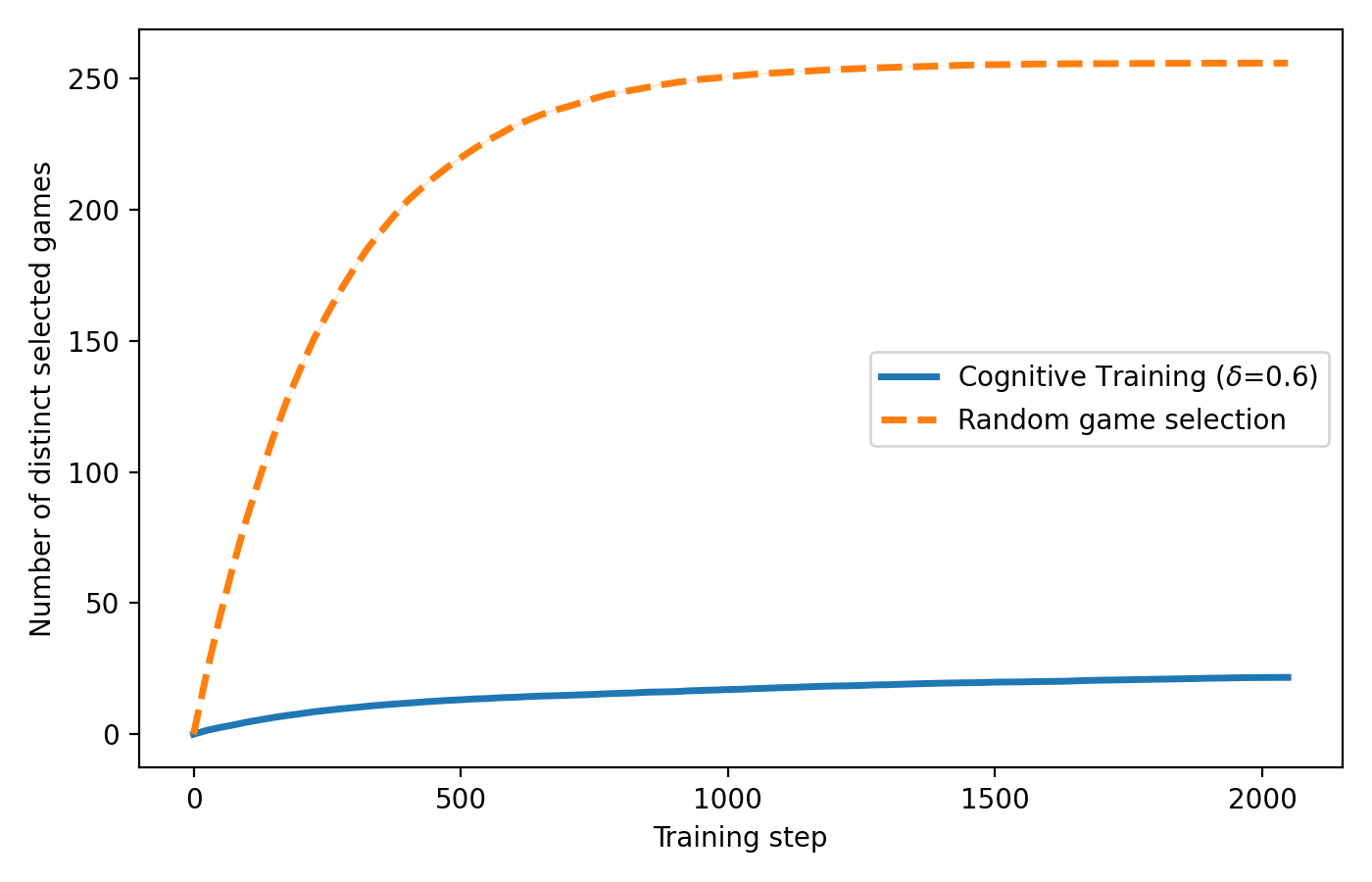}
\caption{Left: Geometric mean of game scores as a function of the diversity
hyper-parameter $\delta$. Middle: geometric mean of game scores during
training. Right: number of distinct games selected during training.
All values are averaged over 64 independently sampled toy worlds;
shaded bands show standard error across worlds.}\label{fig:simulation}
\end{figure}

\subsection{Outlook}

\subsubsection{Playing the Meta-Game}\label{subsec:playing-the-meta-game}

The question of playing\emph{ }the meta-game well, i.e. of how to
train the meta-sampler $\mathcal{M}_{\mathcal{S}}$ to achieve strong
performance on the meta-objective $\mathcal{O}$ is delicate and beyond
the scope of this paper: the space of possible Xent Games is enormous
and enumerating all games is simply out of the question. 

While it is not clear a priori how strong a generalist pre-trained
model needs to be to learn to play $\mathcal{O}$, there is, to the
best of our knowledge, no text game with verifiable rewards that is
playable by a human that an LLM of current sizes (as of Spring 2026)
cannot learn to play at least at a roughly comparable level. A natural
setting to allow for the training is to increase the number of samples.
This can be achieved by e.g. relying on faster transfer learning estimates
to increase sample throughput (such as the upcoming \cite{hongler-emil-renard-gabriel,hongler-asani-gabriel-emil-renard}).
Another promising way is to train small models $\mathcal{M}$ and
use this as a baseline to transfer to larger models.  

\subsubsection{Open-Ended Skill Discovery}\label{subsec:cognitive-training-outlook}

The most significant claim that can be made in the context of cognitive
training is that it can lead to an answer to Problem \ref{prob:find-open-ended-skill-discovery},
i.e. to \emph{relevant skill discovery} (in the sense of Definition
\ref{def:relevant-skill-discovery} above). The claim is based on
the following: 
\begin{itemize}
\item On the absence of examples of verifiable-reward games that humans
can play and that reasonably large language models cannot learn to
play.
\item On the implausibility that we can achieve very high values of $qd$
on Xent Games having $\min\left(q,d\right)$ very small for a nontrivial
range of $\delta$-values $\Delta\subset\left[0,1\right]$.
\end{itemize}
\begin{claim}
\label{claim:relevant-skill-discovery}There is a reasonably large
model $\mathcal{M}$ and meta-sampler model $\mathcal{M}_{\mathcal{S}}$
such that we have the following: for a nontrivial range $\Delta\subset\left[0,1\right]$
of $\delta$-values, cognitive training with external pressure $p=0$
yields unbounded relevant skill discovery. 
\end{claim}

Adding an external benchmark $E$ for which a certain score level
$\lambda_{E}$ is achievable using a curriculum $\mathbb{G}$ (not
necessarily consisting of Xent Games), an analogous question is raised. 
\begin{claim}
\label{claim:relevant-skill-discovery-with-benchmark}Assume $\left(E,\lambda_{E},\mathbb{G}\right)$
is given. There are reasonably large model $\mathcal{M}$ and meta-sampler
model $\mathcal{M}_{\mathcal{S}}$ (with access to $\mathbb{G}$ in
its prompt), a nontrivial range $\Delta\subset\left[0,1\right]$ of
$\delta$-values, and a pressure schedule $p$, such that cognitive
training yields unbounded relevant skill discovery, while converging
to $\lambda_{E}$ performance on $E$. 
\end{claim}

\begin{rem}
The access to $\mathbb{G}$ for $\mathcal{M}_{\mathcal{S}}$ is not
an ideal constraint (which can be voided in natural cases), but it
is needed to avoid some situations where learning $E$ involves access
to e.g. a cryptographic secret. 
\end{rem}

\subsection{Taking a Step Back}\label{subsec:taking-a-step-back}

The derivation of the meta-objective obtained in Section \ref{subsec:meta-objective-derivation}
to get a meta-objective $\mathcal{O}$ for a greedy curriculum is
in fact not very specific to Xent Games. The needed ingredient is
a concrete framework to define training on games written in a specific
language satisfying the following:
\begin{itemize}
\item Composition, i.e. any two games $G,H$ can be concatenated to yield
a game $G\oplus H$ of length $l\left(G\oplus H\right)=l\left(G\right)+l\left(H\right)$.
\item A transfer value definition $\mathcal{T}_{G}^{\mathcal{M}}\left(H\right)$
satisfying: 
\begin{itemize}
\item $\mathcal{T}_{G\oplus G}\left(H\right)\approx2\mathcal{T}_{G}\left(H\right)$
and $\mathcal{T}_{G}\left(H\oplus H\right)=2\mathcal{T}_{G}\left(H\right)$.
\item $\mathcal{T}_{\alpha G}\left(H\right)=\mathcal{T}_{G}\left(H\right)$
for $\alpha>0$ (thus making the explicit support of rescaling unnecessary,
see Remark \ref{rem:important-theoretical-possibility}).
\end{itemize}
\end{itemize}
From this framework, assuming we have a greedy online curriculum of
games that bring relevant new skills to the model, we are naturally
led to a meta-objective of the form
\[
\mathcal{O}=\frac{qd+bp}{l},
\]
where $q,d,b,p$ are as above (see Section \ref{subsec:meta-objective-formula}). 

The derivation's spirit is the following: we \emph{assume a consistent
meta-objective $\mathcal{O}$ exists}, and we find consistency relations
it must satisfy, based on the idea that \emph{the only way in which
a sequence of games matters} is \emph{through its effect on model
training} (and the only way to judge of this is by judging on its
performance on existing games or external metrics). By allowing the
space of games to have enough flexibility, we find that this (somehow
surprisingly) substantially constrains the shape that the meta-objective
must take: if we require that it must not depend on various descriptions
of (effectively) the same curriculum, then we find the above explicit
formula (leaving only two hyper-parameter degrees of freedom). This
line of argument is closely related to the concept of gauge symmetry
in theoretical physics.

\section{Perspectives and Conclusion}

The main claim of this article is informally the following: 
\begin{claim}
If we believe a greedy curriculum of games to build general capabilities
in an LLM can be defined, then it can be replicated by a greedy curriculum
of Xent Games that brings relevant new skills (in the sense of Definition
\ref{def:relevant-skill-discovery}); and if such a curriculum of
Xent Games can be produced, then (under consistency assumptions),
its meta-objective must have the shape of the meta-objective $\mathcal{O}$
given by Expression \ref{eq:meta-objective} in Section \ref{subsec:meta-objective-formula}. 

Cognitive training based on this meta-objective hence seems to offer
a compelling way to build models endowed with general capabilities.
It is worth noting that while the space of Xent Games is a compelling
framework to implement cognitive training, the principles behind it
and the meta-objective derivation can be implemented on any space
of games that is sufficiently rich in terms of compositional structure,
as discussed in Section \ref{subsec:taking-a-step-back} above. 
\end{claim}

\subsection*{Acknowledgements}

Many insights presented in this paper have emerged from collaborations
with Diego Dorn, Vassilis Papadopoulos, Marco Tuccio, Jérémie Wenger,
and Nicolas Zlatoff on Large Language Models, with whom key ideas
were discussed and investigated.

In addition, the authors would like to thank Emmanuel Abbé, Apoorv
Agarwal, Mathieu Alain, Mohammad Asani, Alberto Bietti, Gloria Capano,
Rahul Chalamala, Tarun Chitra, Jordan Cotler, Davide Crapis, Marco
De Rossi, Mario Geiger, Evgenii Golikov, Alex Graves, Nicola Greco,
Bara Hudcová, Arthur Jacot, Niels Linnemann, Tomáš Mikolov, Clément
Moulin-Frier, Max Nye, Mihir Patel, João Penedones, David Pfau, Maciej
Rudzinski, Stanislav Smirnov, Yi Sun, George Walker, Miles Wang, Jason
Wang, Shouqiao Wang, and Matthieu Wyart for interesting discussions,
as well as the participants in the Quine seminar and Demeco workshop
for insightful questions and remarks. 

\section*{Appendix A: S-XGL Specification}

In this appendix, we provide a specification of S-XGL (see Section
\ref{subsec:sxgl}), as implemented on the following repository: \href{https://www.github.com/xentlabs/s-xgl}{https://www.github.com/xentlabs/s-xgl}.
The language specification is deliberately minimal (even more so than
that of O-XGL defined in \cite{hongler-emil}), consisting of 8+1
instructions. A few preliminary remarks are in order: 
\begin{itemize}
\item S-XGL is suitable to reward (and thus train) several models; in the
cognitive training described above, a single model is trained, and
the rewards to other models should just be discarded by the interpreter. 
\item Some design elements can appear cumbersome to a human reader; the
goal is not to make the code particularly human-readable, but to follow
the design elements outlined in Section \ref{subsec:sxgl}.
\item There is no metadata associated with an S-XGL program; S-XGL programs
can be seamlessly concatenated. 
\item Some syntactic sugar can be added to shorten the code without altering
the instruction set; we avoid these questions here. 
\end{itemize}

\subsection*{Global Metadata}

An S-XGL game space $\mathcal{G}$ consists of the following constant
and global (i.e. common to all games $G\in\mathcal{G}$) metadata:
\begin{itemize}
\item a fixed token vocabulary $\mathcal{V}$;
\item a list of models $m_{u}$ for $u<U$, each based on $\mathcal{V}$
(with $m=m_{0}$ being the player model under cognitive training);
\item a list of token strings $s_{k}$ for $k<K$ of common maximal length
$L\geq1$.
\end{itemize}

\paragraph*{Code Structure}
\begin{itemize}
\item It is understood that any S-XGL program ends with an empty line, which
is the clearing instruction, followed by a newline symbol (to enable
seamless concatenation of programs); beyond this constraint, any token
sequence is valid S-XGL code.
\item The interpreter will go line by line through the game code and interpret
any line that matches the syntax of an \emph{instruction line} (i.e.
that satisfies the syntax of one of the 9 instructions below) as described. 
\item Otherwise, the code of the line is not interpreted (though it can
be used to fill token strings by further lines, as described below;
instruction lines can also be used as data to fill token strings).
\end{itemize}

\subsection*{Set of Variables}

Given the global metadata, we define the following set of variables
that are allowed to evolve during a game run:
\begin{itemize}
\item The token strings $s_{k}$ for $k<K$: each token string $s$ consists
$\mathrm{len}\left(s\right)$ tokens in $\mathcal{V}$, with $0\leq\mathrm{len}\left(s\right)\leq l_{\max}$;
at initialization and after clearing, we have $\mathrm{len}\left(s\right)=0$.
\item Each model $m_{u}$ for $u<U$ has a context register and a xent accumulator.
\end{itemize}

\subsection*{The 8+1 possible expressions}
\begin{itemize}
\item $m<<s$: append the string $s$ to the context register of $m$.
\item $s<<m$: elicit $\mathrm{len}\left(s\right)$ tokens from $m$ (given
its current context) and append them to $s$, stopping upon $s$ reaching
the maximal length $l_{\max}$; if $s$ is empty, elicit $1$ token
from $m$. 
\item $s>>m$: compute the xent of $s$ under $m$ (given no context), and
add it to the xent accumulator .
\item $m>>s$: compute the xent of $s$ under $m$ (given no context), and
subtract it from the xent accumulator.
\item $m_{\ell}<<m_{r}$: reward $m_{\ell}$ with the score from the xent
accumulator of $m_{r}$ (the judge) and then clear that xent accumulator.
\item $m_{\ell}>>m_{r}$: reward $m_{r}$ with an ensure nonlinearity (see
Section \ref{subsec:sxg-def-and-basic-runtime-instruction}) applied
to the xent accumulator (the judge) and then clear that xent accumulator. 
\item $s_{\ell}<<s_{r}$: 
\begin{itemize}
\item if $\mathrm{len}\left(s_{r}\right)>0$: append the tokens of $s_{r}$
to those of $s_{\ell}$ and stop upon $s_{l}$ reaching maximal length
$L$. In other words, append the first $\min\left(\mathrm{len}\left(s_{r}\right),l_{\max}-\mathrm{len}\left(s_{\ell}\right)\right)$
tokens from $s_{r}$ to $s_{l}$. 
\item if $\mathrm{len}\left(s_{r}\right)=0$: append the tokenization from
the previous code line $c$ to $s_{l}$ and stop upon $s_{l}$ reaching
maximal length, i.e. append the first $\min\left(\mathrm{len}\left(c\right),l_{\max}-\mathrm{len}\left(s_{r}\right)\right)$
tokens of $c$. 
\end{itemize}
\item $s_{\ell}>>s_{r}$: replace $s_{r}$ with the first $\min\left(\mathrm{len}\left(s_{\ell}\right),\mathrm{len}\left(s_{r}\right)\right)$
tokens of $s_{\ell}$ and remove those tokens from $s_{\ell}$ (shifting
it to the left and reducing its length by $\min\left(\mathrm{len}\left(s_{\ell}\right),\mathrm{len}\left(s_{r}\right)\right)$).
In particular, if $\ell=r$: clear $s_{\ell}$. 
\item empty line: clear all strings, context registers and xent accumulators.
\end{itemize}

\section*{Appendix B: Simulation of Cognitive Training on a Toy Example}

We describe the greedy optimization algorithm used to simulate the
result of Section \ref{subsec:toy-model-for-cognitive-training},
yielding in particular the results shown in Figure \ref{fig:simulation}.
To illustrate the dynamics of the algorithm, and in particular its
behavior with respect to $\delta$, we use a finite-world toy analogue
of the Cognitive Training algorithm. 
\begin{itemize}
\item The game space $\mathcal{G}$ consists of $N_{\mathcal{G}}$ games,
all with the same description length. Each game $G\in\mathcal{G}$
is represented by a non-negative skill vector $(w_{G,s})_{s\in\mathcal{S}}\in\mathbb{R}_{+}^{N_{\mathcal{S}}}$,
where $\mathcal{S}$ is the set of (in practice non-observable) latent
skills. 
\item The game--skill matrix $W=(w_{G,s})$ is sampled as follows: each
entry is sampled from a low-mean exponential distribution $\mathrm{Exp}\left(\mu_{\mathrm{low}}\right)$,
and with low probability $p_{\mathrm{high}}$ we add a sample from
a high-mean exponential distribution $\mathrm{Exp}\left(\mu_{\mathrm{high}}\right)$.
Adding this second component follows the idea that some games strongly
rely on certain skills, and that these skills can have a high transfer
to some other games.
\item A model state is described by a vector of free skills $f=(f_{s})_{s\in\mathcal{S}}$
that is mapped to normalized skills through a saturating transformation
\[
\ensuremath{\mathrm{sk}}_{S}=1-\left(1+\eta f_{s}\right)^{-\alpha_{s}}\in[0,1],
\]
where $\eta>0$ is a scale parameter, and $\alpha_{s}>0$ is a skill-specific
exponent (chosen randomly as $\mu_{\mathrm{skill}}$ times a squared
normal distribution).
\item Given a model state $f,$ the score of the model on the game $G$
is defined by 
\[
\mathrm{S}_{f}(G)=\left(\frac{\sum_{s\in\mathcal{S}}w_{G,s}\ensuremath{\mathrm{sk}}_{s}^{p}}{\sum_{s\in\mathcal{S}}w_{G,s}}\right)^{\frac{1}{p}},
\]
which interpolates between the arithmetic mean ($p=1$) and the weighted
geometric mean $(p\to0)$; in the simulation, we use $p=1$.
\item For a game $G$, training is modeled by the update $f_{s}\mapsto f_{s}+w_{G,s}$,
so that the transfer value from $G$ to $H$ is
\[
\mathcal{T}_{G}^{f}(H)=\mathrm{S}_{f+w_{G}}(H)-\mathrm{S}_{f}(H).
\]
\item We define the performance of the final model as the geometric mean
of the game scores:
\[
\left(\prod_{g\in\mathcal{G}}\mathrm{S}_{f}(g)\right)^{\frac{1}{N_{\mathcal{G}}}}.
\]
\end{itemize}
We do not use external benchmarks and a constant description length
for the games; the meta-objective thus only consists of the internal
relevance term. We apply the Cognitive Training algorithm to construct
the curriculum by greedily selecting, at each step, the game $H$
that maximizes the meta-objective for the chosen value of $\delta.$
In the example illustrated in Figure \ref{fig:simulation}, we run
the experiments with $n_{\mathrm{games}}=256$, $n_{\mathrm{skills}}=64$
$\mu_{\mathrm{low}}=100$, $\mu_{\mathrm{high}}=1$, $p_{\mathrm{high}}^{-1}=\sqrt{n_{\mathrm{games}}n_{\mathrm{skills}}}$,
$\eta=0.01$, $\mu_{\mathrm{skill}}=0.1$ for $n_{\mathrm{steps}}=2048$
steps over $n_{\mathrm{worlds}}=64$ independently sampled toy worlds.


\begin{thebibliography}{XFLZZ25(2025)}
\bibitem[AVZSB24]{agarwal-on-policy-distillation}R. Agarwal, N. Vieillard,
Y. Zhou, P. Stanczyk, S. Ramos, M. Geist, and O. Bachem, On-policy
distillation of language models: Learning from self-generated mistakes.
In The Twelfth International Conference on Learning Representations,
2024. 

\selectlanguage{english}%
\bibitem[AAT2023]{key-1}A.A. Team, J. Bauer, K. Baumli, S. Baveja,
F. Behbahani, A. Bhoopchand, N. Bradley-Schmieg, M. Chang, N. Clay,
A. Collister, V. Dasagi, L. Gonzalez, K. Gregor, E. Hughes, S. Kashem,
M. Loks-Thompson, H. Openshaw, J. Parker-Holder, S. Pathak, N. Perez-Nieves,
N. Rakicevic, T. Rocktäschel, Y. Schroecker, J. Sygnowski, K. Tuyls,
S. York, A. Zacherl, and L. Zhang, Human-Timescale Adaptation in an
Open-Ended Task Space, https://arxiv.org/abs/2301.07608

\bibitem[Bax00]{baxter-inductive-learning}J. Baxter, A Model of Inductive
Bias Learning, \emph{Journal Of Artificial Intelligence Research},
\textbf{12}:149--198, 2000, https://arxiv.org/abs/1106.0245.

\bibitem[BVLFW25]{tutorial-on-meta-rl}J. Beck, R. Vuorio, E.Z. Liu,
Z. Xiong, L. Zintgraf, C. Finn, S. Whiteson, A Tutorial on Meta-Reinforcement
Learning, \emph{Foundations and Trends in Machine Learning} \textbf{18}(2-3):224--384,
https://arxiv.org/pdf/2301.08028.

\bibitem[AFKKQH22]{cont}C. An, J. Feng, K. Lv, L. Kong, X. Qiu, X.
Huang, CoNT: Contrastive Neural Text Generation, \emph{Advances in
Neural Information Processing Systems} \textbf{35} (NeurIPS 2022),
https://arxiv.org/abs/2205.14690.

\bibitem[BaSt07]{banerjee-stone}B. Banerjee and P. Stone, General
Game Learning using Knowledge Transfer, Proceedings of the 20th International
Joint Conference on Artificial Intelligence, 2007. https://www.cs.utexas.edu/\textasciitilde ai-lab/pubs/IJCAI07-bikram.pdf

\bibitem[BCK23]{bluem-czech-kersting}J. Blüm, J. Czech, and K. Kersting,
AlphaZe\textasteriskcentered\textasteriskcentered : AlphaZero-like
baselines for imperfect information games are surprisingly strong,
Front. Artif. Intell., 12 May 2023, Sec. Machine Learning and Artificial
Intelligence, Volume 6, https://doi.org/10.3389/frai.2023.1014561,
2023

\bibitem[CGHCL21]{cghcl-i}S. Carré, F. Gabriel, C. Hongler, G. Lacerda,
and G. Capano, Smart Proofs via Smart Contracts: Succinct and Informative
Mathematical Derivations via Decentralized Markets, https://arxiv.org/abs/2102.03044.

\bibitem[CGHCL24]{cghcl-ii}S. Carré, F. Gabriel, C. Hongler, G. Lacerda,
and G. Capano, Smart Proofs via Recursive Information Gathering: Decentralized
Refereeing by Smart Contracts, \emph{Distributed Ledger Technologies:
Research and Practice}, \textbf{3}(1):1-{}-19 https://dl.acm.org/doi/10.1145/3595298,
2024.

\bibitem[CWWWX23]{survey-llm-eval}Yupeng Chang, Xu Wang, Jindong
Wang, Yuan Wu, Linyi Yang, Kaijie Zhu, Hao Chen, Xiaoyuan Yi, Cunxiang
Wang, Yidong Wang, Wei Ye, Yue Zhang, Yi Chang, Philip S. Yu, Qiang
Yang, Xing Xie, A Survey on Evaluation of Large Language Models, https://arxiv.org/abs/2307.03109

\bibitem[CZJGS24]{chatbot-arena}W.-L. Chiang, L. Zheng, Y. Sheng,
A.N. Angelopoulos, T. Li, D. Li, H. Zhang, B. Zhu, M. Jordan, J.E.
Gonzalez, I. Stoica, Chatbot Arena: An Open Platform for Evaluating
LLMs by Human Preference, https://arxiv.org/abs/2403.04132.

\bibitem[Cho19]{chollet}F. Chollet, On the Measure of Intelligence,
https://arxiv.org/abs/1911.01547.

\bibitem[CKKLP25]{chollet-knoop-kamradt-landers-pinkard}F. Chollet,
M. Knoop, G. Kamradt, B. Landers, H. Pinkard, ARC-AGI-2: A New Challenge
for Frontier AI Reasoning Systems, https://arxiv.org/abs/2505.11831

\bibitem[Clu19]{clune}J. Clune, AI-GAs: AI-generating algorithms,
an alternate paradigm for producing general artificial intelligence,
https://arxiv.org/abs/1905.10985v2.

\bibitem[CKATT18]{textworld}M.-A. Côté, A. Kádár, X. Yuan, B. Kybartas,
T. Barnes, E. Fine, J. Moore, R.Y. Tao, M. Hausknecht, L.E. Asri,
M. Adada, W. Tay, and A. Trischler, TextWorld: A Learning Environment
for Text-based Games, Proceedings of the Computer Games Workshop,
International Joint Conference on Artificial Intelligence 2018, https://arxiv.org/abs/1806.11532

\selectlanguage{american}%
\bibitem[CoTh06]{cover-thomas}\foreignlanguage{english}{T.M. Cover,
J.A. Thomas, Elements of Information Theory (2nd Edition), Wiley,
2006}

\selectlanguage{english}%
\bibitem[DAW24]{das-amini-wu}B.C. Das, M. H. Amini, and Y. Wu, Security
and Privacy Challenges of Large Language Models: A Survey, https://arxiv.org/abs/2402.00888.

\selectlanguage{american}%
\bibitem[Dee24]{deepseek-v3}\foreignlanguage{english}{DeepSeek-V3
Technical Report, DeepSeek-AI, https://arxiv.org/abs/2412.19437.}

\bibitem[DDSW25]{reinforcement-pretraining}\foreignlanguage{english}{Q.
Dong, L. Dong, Y. Tang, T. Ye, Y. Sun, Z. Sui, F. Wei, Reinforcement
Pre-Training, https://arxiv.org/pdf/2506.08007}

\selectlanguage{english}%
\bibitem[ECMO23]{etcheverry-chan-moulin-frier-oudeyer}M. Etcheverry,
B. W.-C. Chan, C. Moulin-Frier, P.-Y. Oudeyer, Meta-Diversity Search
in Complex Systems, A Recipe for Artificial Open-Endedness? https://arxiv.org/abs/2312.00455

\bibitem[GHWZ16]{ghani-hedges-winschel-zahn}N. Ghani, J. Hedges,
V. Winschel, P. Zahn, Compositional game theory, Proceedings of the
33rd Annual ACM/IEEE Symposium on Logic in Computer Science. LICS
'18. New York, NY, US: ACM. pp. 472--481. https://arxiv.org/abs/1603.04641

\bibitem[GGKPW25]{gollakota-gopala-kara-peale-wieder}A. Gollakota,
P. Gopalan, A. Karan, P. Peale, U. Wieder, When does a predictor know
its own loss?, https://arxiv.org/abs/2502.20375.

\bibitem[GBMMK17]{automated-curriculum}A. Graves, M.G. Bellemare,
J. Menick, R. Munos, K. Kavukcuoglu, Automated Curriculum Learning
for Neural Networks, Proceedings of the 34th International Conference
on Machine Learning, Sydney, Australia, PMLR 70, 2017. https://arxiv.org/abs/1704.03003. 

\selectlanguage{american}%
\bibitem[HAPCC25]{hatamizadeh-akter-prabhumoy-kautz-et-al}\foreignlanguage{english}{Ali
Hatamizadeh, Syeda Nahida Akter, Shrimai Prabhumoye, Jan Kautz, M.
Patwary, M. Shoeybi, B. Catanzaro, Y. Choi, RLP: Reinforcement as
a Pretraining Objective, https://arxiv.org/abs/2510.01265.}

\selectlanguage{english}%
\bibitem[HiVC93]{hinton-van-camp} G. E. Hinton and D. Van Camp, Keeping
the neural networks simple by minimizing the description length of
the weights. In Proceedings of the sixth annual conference on Computational
learning theory, pp. 5--13. ACM, 1993. https://dl.acm.org/doi/10.1145/168304.168306

\bibitem[HoEm25]{hongler-emil}C. Hongler and A. Emil, Cross-Entropy
Games for Language Models: From Implicit Knowledge to General Capability
Measures, https://arxiv.org/abs/2506.06832.

\selectlanguage{american}%
\bibitem[HERG26]{hongler-emil-renard-gabriel}C. Hongler, A. Emil,
A. Renard, F. Gabriel, Frost Scores and Cross-Entropy Games: Sample-Efficient
Training, in preparation. 

\selectlanguage{english}%
\bibitem[HAGER26]{hongler-asani-gabriel-emil-renard}C. Hongler, A.
Renard, F. Gabriel, A. Emil, Higher-Order Frost Scoring Schemes, in
preparation.

\selectlanguage{american}%
\bibitem[HLBBG26]{hubotter-et-al26-self-distillation} J. Hübotter,
F. Lübeck, L. Behric, A. Baumann, M. Bagatella, D. Marta, I. Hakimi,
I. Shenfeld, T. K. Buening, C. Guestrin, A. Krause. Reinforcement
Learning via Self-Distillation, https://arxiv.org/pdf/2601.20802. 

\bibitem[HWMHS23]{token-level-adversarial-prompt-detection}Z. Hu,
G. Wu, S. Mitra, R. Zhang, T. Sun, H. Huang, V. Swaminathan, Token-Level
Adversarial Prompt Detection Based on Perplexity Measures and Contextual
Information, https://arxiv.org/abs/2311.11509.

\selectlanguage{english}%
\bibitem[HDPHR24]{hughes-dennis-parker-holder-bekhabani-mavalankar-shi-schaul-rocktaechel}E.
Hughes, M.D. Dennis, J. Parker-Holder, F. Behbahani, A. Mavalankar,
Y. Shi, T. Schaul, T. Rocktäschel, Open-Endedness is Essential for
Artificial Superhuman Intelligence, \emph{Proceedings of the 41st
International Conference on Machine Learning}, PMLR \textbf{235}:20597-20616,
2024. https://arxiv.org/abs/2406.04268

\selectlanguage{american}%
\bibitem[Hut05]{hutter}\foreignlanguage{english}{M. Hutter, \emph{Universal
algorithmic intelligence: Sequential Decisions based on Algorithmic
Probability}. Springer, Berlin, 2005.}

\bibitem[JTZA26]{Tutunov-Scalable-Power-Sampling-26}\foreignlanguage{english}{X.
Ji, R. Tutunov, M. Zimmer, H. B. Ammar. Scalable Power Sampling: Unlocking
Efficient, Training-Free Reasoning for LLMs via Distribution Sharpening.
https://www.arxiv.org/pdf/2601.21590. }

\bibitem[KCAHK22]{Language Models (Mostly) Know What They Know}S.\foreignlanguage{english}{
Kadavath, T. Conerly, A. Askell, T. Henighan, D. Drain, E. Perez,
N. Schiefer, Z. Hatfield-Dodds, N. DasSarma, E. Tran-Johnson, S. Johnston,
S. El-Showk, A. Jones, N. Elhage, T. Hume, A. Chen, Y. Bai, S. Bowman,
S. Fort, D. Ganguli, D. Hernandez, J. Jacobson, J. Kernion, S. Kravec,
L. Lovitt, K. Ndousse, C. Olsson, S. Ringer, D. Amodei, T. Brown,
J. Clark, N. Joseph, B. Mann, S. McCandlish, C. Olah, and J. Kaplan,
Language Models (Mostly) Know What They Know, https://arxiv.org/pdf/2207.05221,
2022}

\bibitem[KMHBA20]{scaling-laws}J.\foreignlanguage{english}{ Kaplan,
S. McCandlish, T. Henighan, T. B. Brown, B. Chess, R. Child, S. Gray,
A. Radford, J. Wu, D. Amodei, Scaling Laws for Neural Language Models,
https://arxiv.org/abs/2001.08361.}

\selectlanguage{english}%
\bibitem[KGRMI22]{kojima-gu-reid-matsuo-iwasawa}T. Kojima, S.S. Gu,
M. Reid, Y. Matsuo, Y. Iwasawa, Large Language Models are Zero-Shot
Reasoners, https://arxiv.org/abs/2205.11916

\bibitem[Lal2015]{laland-et-al}K.N. Laland, T. Uller, M.W. Feldman,
K. Sterelny, B.B. Müller, A. Moczek, E. Jablonka, and J. Odling-Smee,
The extended evolutionary synthesis: its structure, assumptions and
predictions, \emph{Proc. R. Soc. B}, \textbf{282}: 20151019, 2015.

\selectlanguage{american}%
\bibitem[LeHu07]{legg-hutter}\foreignlanguage{english}{S. Legg, M.
Hutter, Universal Intelligence: A Definition of Machine Intelligence,
\emph{Minds \& Machines}, 17(4):391-{}-444, https://arxiv.org/abs/0712.3329,
2007. }

\bibitem[LeHu06]{legg-hutter-2}\foreignlanguage{english}{S. Legg,
M. Hutter, A Formal Measure of Machine Intelligence, IDSIA Technical
Report 10-06, arXiv:cs/0605024v1}

\selectlanguage{english}%
\bibitem[LeSt10]{lehman-stanley-efficiently-evolving}J. Lehman and
K.O. Stanley, Efficiently evolving programs through the search for
novelty, \emph{GECCO '10: Proceedings of the 12th annual conference
on Genetic and evolutionary computation}, 2010, https://doi.org/10.1145/1830483.1830638. 

\bibitem[LeSt11a]{lehmann-stanley-evolution-novelty}J. Lehman and
K.O. Stanley. Abandoning objectives: Evolution through the search
for novelty alone. Evolutionary Computation, 19(2):189--223, 2011.

\bibitem[LeSt11b]{lehmann-stanley-novelty-search}J. Lehman and K.O.
Stanley. Evolving a diversity of virtual creatures through novelty
search and local competition. In Proceedings of the 13th annual conference
on Genetic and evolutionary computation, pages 211--218. ACM, 2011.

\bibitem[LBL23]{holistic}P. Liang, R. Bommasani, T. Lee, D. Tsipras,
D. Soylu, M. Yasunaga, Y. Zhang, D. Narayanan, Y. Wu, A. Kumar, B.
Newman, B. Yuan, B. Yan, C. Zhang, C. Cosgrove, C. D. Manning, C.
Ré, D. Acosta-Navas, D. A. Hudson, E. Zelikman, E. Durmus, F. Ladhak,
F. Rong, H. Ren, H. Yao, J. Wang, K. Santhanam, L. Orr, L. Zheng,
M. Yuksekgonul, M. Suzgun, N. Kim, N. Guha, N. Chatterji, O. Khattab,
P. Henderson, Q. Huang, R. Chi, S. M. Xie, S. Santurkar, S. Ganguli,
T. Hashimoto, T. Icard, T. Zhang, V. Chaudhary, W. Wang, X. Li, Y.
Mai, Y. Zhang, Y. Koreeda, Holistic Evaluation of Language Models,
https://arxiv.org/abs/2211.09110

\bibitem[LHFLL22]{contrastive}X.L. Li, A. Holtzman, D. Fried, P.
Liang, J. Eisner, T. Hashimoto, L. Zettlemoyer, M. Lewis, Contrastive
Decoding: Open-ended Text Generation as Optimization, Association
for Computer Linguistics, https://arxiv.org/abs/2210.15097

\selectlanguage{american}%
\bibitem[LKKYW25]{spice-meta}B. Liu, C. Jin, S. Kim, W. Yuan., W.
Zhao, I. Kulikov, X. Li, S. Sukhbaatar, J. Lanchantin, J. Weston,
SPICE: Self-Play In Corpus Environments Improves Reasoning. https://arxiv.org/pdf/2510.24684,
2025. 

\bibitem[LJJSG25]{data-sentinel}Y.\foreignlanguage{english}{ Liu,
Y. Jia, J. Jia, D. Song, N.Z. Gong, DataSentinel: A Game-Theoretic
Detection of Prompt Injection Attacks, 2025 IEEE Symposium on Security
and Privacy (SP), 2190-2208, 2025.}

\bibitem[LYZZX24]{automatic-universal-prompt-injection-attacks}X.\foreignlanguage{english}{
Liu},\foreignlanguage{english}{ }Z.\foreignlanguage{english}{ Yu},\foreignlanguage{english}{
}Y.\foreignlanguage{english}{ Zhang},\foreignlanguage{english}{ }N.\foreignlanguage{english}{
Zhang},\foreignlanguage{english}{ }C.\foreignlanguage{english}{ Xiao,
Automatic and Universal Prompt Injection Attacks against Large Language
Models, https://arxiv.org/abs/2403.04957.}

\selectlanguage{english}%
\bibitem[LWLZD24]{best-practices-synth-data}R. Liu, J. Wei, F. Liu,
C. Si, Y. Zhang, J. Rao, S. Zheng, D. Peng, D. Yang, D. Zhou, A. M.
Dai, Best Practices and Lessons Learned on Synthetic Data, COLM 2024,
https://arxiv.org/abs/2404.07503

\selectlanguage{american}%
\bibitem[LTML25]{thinking-machines-on-policy}\foreignlanguage{english}{K.
Lu and Thinking Machines Lab, \textquotedbl On-Policy Distillation\textquotedbl ,
Thinking Machines Lab: Connectionism, Oct 2025.}

\bibitem[MPA24]{mbuya-pfosser-anastosopoulos}\foreignlanguage{english}{J.
K. Mbuya, D. Pfoser, A. Anastasopoulos, Trajectory Anomaly Detection
with Language Models, \emph{SIGSPATIAL '24: Proceedings of the 32nd
ACM International Conference on Advances in Geographic Information
Systems}, 208--219, https://arxiv.org/abs/2409.15366. }

\selectlanguage{english}%
\bibitem[MCV20]{meister-cotterell-vieira}Clara Meister, Ryan Cotterell,
and Tim Vieira. If beam search is the answer, what was the question?
\emph{Proceedings of the 2020 Conference on Empirical Methods in Natural
Language Processing (EMNLP)}, pp. 2173--2185. Association for Computational
Linguistics, 2020. doi: 10.18653/v1/2020.emnlp-main.170. https://aclanthology.org/2020.emnlp-main.170

\bibitem[MKSLH15]{atari}V. Mnih, K. Kavukcuoglu, D. Silver, A.A.
Rusu, J. Veness, M.G. Bellemare, A. Graves, M. Riedmiller, A.K. Fidjeland,
G. Ostrovski, S. Petersen, C. Beattie, A. Sadik, I. Antonoglou, H.
King, D. Kumaran, D. Wierstra, S. Legg, D. Hassabis, Human-level control
through deep reinforcement learning, Nature 518, 529--533 (2015).
https://doi.org/10.1038/nature14236

\selectlanguage{american}%
\bibitem[MJB15]{mikolov-joulin-baroni}T. Mikolov, A. Joulin, M. Baroni,
A Roadmap towards Machine Intelligence, \foreignlanguage{english}{https://arxiv.org/abs/}1511.08130v2.

\selectlanguage{english}%
\bibitem[MSBLB17]{deepstack}M. Moravčík, M. Schmid, N. Burch, V.
Lisý, D. Morrill, N. Bard, T. Davis, K. Waugh, M. Johanson, M. Bowling,
DeepStack: Expert-Level Artificial Intelligence in No-Limit Poker,
https://arxiv.org/abs/1701.01724. 

\bibitem[OAI2020]{openai-2020}T.B. Brown, B. Mann, N. Ryder, M. Subbiah,
J. Kaplan, P. Dhariwal, A. Neelakantan, P. Shyam, G. Sastry, A. Askell,
S. Agarwal, A. Herbert-Voss, G. Krueger, T. Henighan, R. Child, A.
Ramesh, D.M. Ziegler, J. Wu, C. Winter, C. Hesse, M. Chen, E. Sigler,
M. Litwin, S. Gray, B. Chess, J. Clark, C. Berner, S. McCandlish,
A. Radford, I. Sutskever, D. Amodei, Language Models are Few-Shot
Learners, \emph{Advances in Neural Information Processing Systems}
\textbf{33}:1877--1901, https://arxiv.org/abs/2005.14165.

\bibitem[PBRW99]{page-brin-motwani-winograd}L. Page, S. Brin, R.
Motwani, and T. Winograd, The PageRank Citation Ranking: Bringing
Order to the Web. \emph{Technical Report. Stanford InfoLab.} http://ilpubs.stanford.edu:8090/422/

\bibitem[PWH24]{papadopoulos-wenger-hongler}V. Papadopoulos, J. Wenger,
C. Hongler, Arrows of Time for Large Language Models, International
Conference on Machine Learning 2024. https://arxiv.org/abs/2401.17505,
2024.

\bibitem[PBS16]{parisotto-ba-salakhutidnov}E. Parisotto, J.L. Ba,
R. Salakhutdinov, Actor-Mimic: Deep Multitask and Transfer Reinforcement
Learning, International Conference on Learning Representations, 2016,
https://arxiv.org/abs/1511.06342.

\selectlanguage{american}%
\bibitem[PSS16]{quality-diversity}J.\foreignlanguage{english}{K.
Pugh, L.B. Soros, K. O. Stanley, Quality Diversity: A New Frontier
for Evolutionary Computation, \emph{Front. Robot. AI,} 3 - 2016 |
https://doi.org/10.3389/frobt.2016.00040.}

\selectlanguage{english}%
\bibitem[PCDSC24]{useful-llm-eval}Ji-Lun Peng, Sijia Cheng, Egil
Diau, Yung-Yu Shih, Po-Heng Chen, Yen-Ting Lin, Yun-Nung Chen, A Survey
of Useful LLM Evaluation, https://arxiv.org/abs/2406.00936v1

\bibitem[Schmi07]{schmidhuber}J. Schmidhuber. Gödel machines: Fully
self-referential optimal universal self-improvers. \emph{ Artificial
general intelligence}, pages 199--226. Springer, 2007. https://sferics.idsia.ch/pub/juergen/gmAGI.pdf

\bibitem[Schmi10]{schmidhuber-survey}J. Schmidhuber. Formal Theory
of Creativity, Fun, and Intrinsic Motivation (1990-2010). \emph{IEEE
Transactions on Autonomous Mental Development}, 2(3):230-247, 2010.
IEEE, https://people.idsia.ch/\textasciitilde juergen/ieeecreative.pdf

\bibitem[SYCYL24]{shi-yang-cai-zhang-wang-yang-lam}Chufan Shi, Haoran
Yang, Deng Cai, Zhisong Zhang, Yifan Wang, Yujiu Yang, Wai Lam, A
Thorough Examination of Decoding Methods in the Era of LLMs, https://arxiv.org/abs/2402.06925

\bibitem[SYSKZ24]{meta-learning-llm}Sanchit Sinha, Yuguang Yue, Victor
Soto, Mayank Kulkarni, Jianhua Lu, Aidong Zhang, MAML-en-LLM: Model
Agnostic Meta-Training of LLMs for Improved In-Context Learning, International
Conference on Knowledge Discovery and Data Mining, 2024, https://arxiv.org/abs/2405.11446

\bibitem[Sig23]{sigaud-baldassare-colas-doncieux-duro-oudeyer-perrin-gilbert-santucci}O.
Sigaud, G. Baldassarre, C. Colas, S. Doncieux, R. Duro, P.-Y. Oudeyer,
N. Perrin-Gilbert, V.G. Santucci, A Definition of Open-Ended Learning
Problems for Goal-Conditioned Agents, https://arxiv.org/abs/2311.00344.

\bibitem[SHSH18]{alpha-zero}D. Silver, T. Hubert, J. Schrittwieser,
I. Antonoglou, M. Lai, A. Guez, M. Lanctot, L. Sifre, D. Kumaran,
T. Graepel, T. Lillicrap, K. Simonyan, and D. Hassabis. A general
reinforcement learning algorithm that masters chess, shogi, and Go
through self-play, \emph{Science}, 7 Dec 2018, \textbf{362}(6419):1140-1144,
DOI: 10.1126/science.aar6404, 2018.

\selectlanguage{american}%
\bibitem[SDHA26]{shenfeld-self-distillation}I. Shenfeld, M. Damani,
J. Hübotter, P. Agrawal. Self-Distillation Enables Continual Learning.
https://arxiv.org/pdf/2601.19897. 

\selectlanguage{english}%
\bibitem[SWLL24]{song-wang-li-lin}Y. Song, G. Wang, S. Li, B.Y. Lin,
The Good, The Bad, and The Greedy: Evaluation of LLMs Should Not Ignore
Non-Determinism, https://arxiv.org/abs/2407.10457. 

\bibitem[StBy19]{stahlberg-byrne}Felix Stahlberg and Bill Byrne.
On NMT search errors and model errors: Cat got your tongue? Proceedings
of the 2019 Conference on Empirical Methods in Natural Language Processing
and the 9th International Joint Conference on Natural Language Processing
(EMNLP-IJCNLP), pp. 3356--3362. Association for Computational Linguistics,
2019. doi: 10.18653/v1/D19-1331. URL https://aclanthology.org/D19-1331.

\bibitem[SSSJ25]{stojanovski-stanley-sharratt-jones-adefioye-kaddour-koepf}Z.
Stojanovski, O. Stanley, J. Sharratt, R. Jones, A. Adefioye, Jean
Kaddour, Andreas Köpf, REASONING GYM: Reasoning Environments for Reinforcement
Learning with Verifiable Rewards, https://arxiv.org/abs/2505.24760 

\bibitem[TKM23]{tan-kazemi-mihalcea}Q. Tan, A. Kazemi, R. Mihalcea,
Text-Based Games as a Challenging Benchmark for Large Language Models,
International Conference on Learning Representations, TinyPapers,
2023. 

\bibitem[Tur50]{turing}A. M. Turing, Computing Machinery and Intelligence.
\emph{Mind} \textbf{49}:433-460, 1950.

\selectlanguage{american}%
\bibitem[WLCS19]{POET}\foreignlanguage{english}{R. Wang, J. Lehman,
J. Clune, K. O. Stanley. POET: open-ended coevolution of environments
and their optimized solutions. In \emph{Proceedings of the Genetic
and Evolutionary Computation Conference}, pp. 142--151, ACM, 2019. }

\selectlanguage{english}%
\bibitem[WWZ24]{data-augmentation}Z. Wang, P. Wang, K. Liu, P. Wang,
Y. Fu, C.-T. Lu, C.C. Aggarwal, J. Pei, Y. Zhou, A Comprehensive Survey
on Data Augmentation, https://arxiv.org/abs/2405.09591

\selectlanguage{american}%
\bibitem[WXSLD22]{cot}\foreignlanguage{english}{J. Wei, X. Wang,
D. Schuurmans, M. Bosma, B. Ichter, F. Xia, E. Chi, Q. Le, D. Zhou,
Chain-of-Thought Prompting Elicits Reasoning in Large Language Models,
Proceedings of NeurIPS 2022, https://arxiv.org/abs/2201.11903.}

\selectlanguage{english}%
\bibitem[Wil92]{williams}R.J. Williams, Simple Statistical Gradient-Following
Algorithms for Connectionist Reinforcement Learning, \emph{Machine
Learning}, \textbf{8}(3-4):229 - 256, https://dl.acm.org/doi/10.1007/bf00992696.

\bibitem[XDCGZ24]{agent-gym}Z. Xi, Y. Ding, W. Chen, B. Hong, H.
Guo, J. Wang, D. Yang, C. Liao, X. Guo, W. He, S. Gao, L. Chen, R.
Zheng, Y. Zou, T. Gui, Q. Zhang, X. Qiu, X. Huang, Z. Wu, Y.-G. Jiang,
AgentGym: Evolving Large Language Model-based Agents across Diverse
Environments, https://arxiv.org/abs/2406.04151v1.

\bibitem[XFLZZ25(2025)]{xing-fan-lou-li-zhang-zhang}X. Xing, Z. Fan,
J. Lou, G. Li, J. Zhang, D. Zhang, PretrainZero: Reinforcement Active
Pretraining, https://arxiv.org/pdf/2512.03442

\bibitem[ZHLLC25]{zhang-hu-lu-lange-clune}J. Zhang, S. Hu, C. Lu,
R. Lange, J. Clune, Darwin Gödel Machine: Open-Ended Evolution of
Self-Improving Agents, https://arxiv.org/abs/2505.22954.\selectlanguage{american}%

\end{thebibliography}
\end{document}